\documentclass[a4paper,12pt,leqno]{amsart}

\usepackage{amsmath,amssymb,amsthm,latexsym,mathrsfs,times}

\newtheorem{theorem}{Theorem}[section]
\newtheorem{lemma}[theorem]{Lemma}
\newtheorem{proposition}[theorem]{Proposition}
\newtheorem{corollary}[theorem]{Corollary}
\newtheorem{definition}{Definition}

\numberwithin{equation}{section}
\allowdisplaybreaks[4]

\begin{document}

\title[Microlocal smoothing effect in a Gevrey class]
{Microlocal smoothing effect \\
for the Schr\"odinger evolution 
equation \\
in a Gevrey class}
\author[R. Mizuhara]{Ryuichiro Mizuhara}
\address{Mathematical Institute, Tohoku University, Sendai 980-8578,
Japan}
\keywords{Microlocal smoothing effect, 
Gevrey class, Almost analytic extension. }
\email{s98m31@math.tohoku.ac.jp}
\subjclass[2000]{Primary 35B65, Secondary 35Q40}

\begin{abstract}
We discuss the microlocal Gevrey smoothing effect 
for the Schr\"odinger equation with variable coefficients
via the propagation property of the wave front set of 
homogenous type. 
We apply the microlocal exponential estimates in a Gevrey case 
to prove our result.   
\end{abstract}

\maketitle

\section{Introduction}\label{sec:intro}

In this paper we study the microlocal smoothness 
in a Gevrey class of solutions 
to the time-dependent Schr\"odinger equation 
with variable coefficients. 

Let $P$ be a Schr\"odinger operator in $\mathbb R^n$ 
\[
P=\frac{1}{2}\sum_{j,k=1}^n a_{jk}(x)D_j D_k
+\sum_{j=1}^n b_j(x)D_j+c(x)
\quad \left( D_j=-i\frac{\partial}{\partial x_j}\right).
\]
We assume that the coefficients of $P$ satisfy the following conditions.

\vspace*{12pt}

\noindent
{\bf Assumption (A)}\enskip
We assume that 
\begin{itemize}
\item 
$a_{jk}(x)\in C^{\infty}(\mathbb R^n;\mathbb R)
\quad(1 \leq j,k \leq n)$,
\item
$b_j(x) \in C^{\infty}(\mathbb R^n)
\quad (1 \leq j \leq n)$, 
\quad 
$c(x) \in C^{\infty}(\mathbb R^n)$.
\item 
The matrix $(a_{jk}(x))_{1 \le j,k \le n}$ 
is symmetric and positive 
definite.
\item
There exist $s>1$, $\sigma>0$, $C_0>0$ and $K_0>0$ such that
\begin{align*}
|\partial_x^{\alpha}(a_{jk}(x)-\delta_{jk})|
& \le C_0 K_0^{|\alpha|}\alpha!^s 
\langle x \rangle^{-\sigma-|\alpha|}
\quad (1 \le j,k \le n),\\ 
|\partial_x^{\alpha} \, \mathrm{Re} \, b_j (x)|
& \le C_0 K_0^{|\alpha|}\alpha!^s 
\langle x \rangle^{1-\sigma-|\alpha|}
\quad (1 \le j \le n),\\
|\partial_x^{\alpha} \, \mathrm{Im} \, b_j (x)|
& \le C_0 K_0^{|\alpha|}\alpha!^s 
\langle x \rangle^{1/s-1-\sigma-|\alpha|}
\quad (1 \le j \le n),\\ 
|\partial_x^{\alpha} \, \mathrm{Re} \, c(x)|
& \le C_0 K_0^{|\alpha|}\alpha!^s 
\langle x \rangle^{2-\sigma-|\alpha|}, \\ 
|\partial_x^{\alpha} \, \mathrm{Im} \, c(x)|
& \le C_0 K_0^{|\alpha|}\alpha!^s 
\langle x \rangle^{1/s-\sigma-|\alpha|}
\end{align*} 
for $\alpha \in \mathbb Z_{+}^n=(\mathbb N \cup \{0\})^n$, 
$x \in \mathbb R^n$, 
where 
$\mathbb N=\{1,2,3, \cdots \}$, 
$\langle x \rangle =(1+|x|^2)^{1/2}$ and 
$\delta_{jk}$ is the Kronecker's delta.
\end{itemize}

\vspace*{12pt}

\noindent{\bf Remark.}\enskip
We can assume $0<\sigma \le 1/s$ without loss of generality.

\vspace*{12pt}

Let $T>0$. 
We consider the solution 
$u(t,\cdot)\in C([0,T];L^2(\mathbb R^n))$ 
to the initial value problem
\begin{equation}\label{eqn:ivp}
\begin{cases}
\dfrac{\partial u}{\partial t}+iPu=0 \quad 
&\text{in} \quad (0,T)\times \mathbb R^n, \\
u(0,x)=u_0(x) \quad &\text{in} \quad \mathbb R^n,  
\end{cases} 
\end{equation}
where $u_0 \in L^2(\mathbb R^n).$

It is well-known that solutions to the Schr\"odinger
equation satisfy the property called a smoothing effect
or gain of regularity, that is, 
the decay of the initial data implies 
the regularity of the solution at positive time. 
Furthermore, more the data decays 
more the solution is smooth.

After the pioneering work of Craig-Kappeler-Strauss
\cite{CKS}, 
the microlocal structure of this phenomenon 
in the variable coefficients case, 
or more precisely the relation between 
the microlocal regularity of solution 
and the behavior of the initial data 
along the backward bicharacteristic, 
studied extensively: see e.g.,
\cite{Do00}, 
\cite{It}, \cite{Na05}, \cite{Wu}
for the $C^{\infty}$ case, 
\cite{MNS06},\cite{MRZ},\cite{RZ99},\cite{RZ00},\cite{RZ02} 
for the analytic case, and 
\cite{KT} for the Gevrey case.

In particular, Nakamura \cite{Na05} introduced 
a new notion of wave front set, 
the homogenous wave front set, 
and extended the result of
Craig-Kappeler-Strauss \cite{CKS} 
to the case of long-range perturbations
by employing the propagation theorem of homogenous 
wave front set. 
The homogenous wave front set
is a conic set in the phase space. 
It propagates along free classical trajectories and 
is suitable to describe the singularity of a solution to 
free Schr\"odinger equations. 

It is remarked that Ito \cite{It} showed 
the essentially equivalence
of the homogenous wave front set and
the quadratic scattering wave front set,
which is a notion for problems on scattering manifolds,
used by Wunsch \cite{Wu} and Hassell-Wunsch \cite{HW}.

For the analytic category problem
with asymptotically flat metrics, 
Martinez-Nakamura-Sordoni \cite{MNS06}
introduced analytic homogenous wave front set 
and generalized the 
results of Robbiano-Zuily \cite{RZ99}, \cite{RZ00}
by simple proofs.

We also remark that 
the results on the characterization
of the wave front set are obtained by   
Hassell-Wunsch \cite{HW}, Ito-Nakamura\cite{IN}, 
Martinez-Nakamura-Sordoni \cite{MNS07} and Nakamura \cite{Na06}.

The purpose of this paper is to refine the microlocal
Gevrey smoothing phenomenon for the Schrodinger equation 
from a view point of the propagation of 
Gevrey wave front set of homogenous type.
We shall prove that a theorem similar to 
the $C^{\infty}$ case or the analytic case 
holds for the Gevrey case.  
Following \cite {MNS06}, we employ microlocal energy method
to prove our result. 
More precisely, combining almost analytic extension of 
Gevrey symbols due to Jung \cite{Ju} 
and microlocal exponential weighted estimates, 
we show the exponential decay of solutions in some direction
on the phase space
under the appropriate condition of the initial data.

We remark that our result is not a generalization of 
the work of Kajitani-Taglialatela \cite{KT}
in which they employ the Fourier integral operator 
with complex-valued phase function.
Indeed, they did not assume
the asymptotically flatness for the metric.
However we emphasize that 
our result can apply the case with the unbounded or 
complex-valued lower order terms.

We here introduce the two kinds of wave front set. 
Let $h,\mu>0$. 
Let $T_{h,\mu}: L^2(\mathbb R^n)\to L^2(\mathbb R^{2n})$ be 
the global FBI transform (the Bargmann transform)
defined by 
\begin{equation}\label{eqn:fbi}
T_{h,\mu}u(x,\xi)
=c_{h,\mu}\int_{\mathbb R^n}
e^{i(x-y)\cdot \xi/h-\mu|x-y|^2/2h}u(y)\,dy
\end{equation}
$(c_{h,\mu}=2^{-n/2}\mu^{n/4}(\pi h)^{-3n/4}).$
\begin{definition}\label{def:wf}
Let $s>1$, $\mu>0$, 
$(x_0, \xi_0)\in \mathbb R^n \times 
(\mathbb R^n \setminus \{0\})$, 
and 
$u \in \mathscr{S}'(\mathbb R^n).$
The point $(x_0,\xi_0)$ 
does not belong to 
the Gevrey wave front set of order $s$ of $u$ 
$((x_0,\xi_0)\notin \mathrm{WF}_s(u))$ 
if there exist a neighborhood $U$ of $(x_0,\xi_0)$
and $\delta>0$, $C>0$ such that 
\[
 \|T_{h,\mu}u\|_{L^2(U)}\le C \exp(-\delta/h^{1/s})
\quad (0 <h \le 1).
\]
\end{definition}
\begin{definition}\label{def:hwf}
Let $s>1$, $\mu>0$, 
$(x_0, \xi_0)\in \mathbb R^{2n}\setminus\{0\}$
and
$u \in \mathscr{S}'(\mathbb R^n).$
The point $(x_0,\xi_0)$ 
does not belong to 
the Gevrey homogenous wave front set of order $s$ 
of $u$ 
$((x_0,\xi_0)\notin \mathrm{HWF}_s(u))$ 
if there exist a conic neighborhood $\Sigma$ 
of $(x_0,\xi_0)$
and $\delta>0$ such that 
\[
 \|\exp\{\delta(|x|^{1/s}+|\xi|^{1/s})\}
T_{1, \mu} u\|_{L^2(\Sigma)}
<+\infty.
\]
\end{definition}
\noindent{\bf Remark.}\enskip
In the case $s=1$, they coincide with 
the analytic wave front set $\mathrm{WF}_a$
and the analytic homogenous wave front set $\mathrm{HWF}_a$
introduced in \cite{MNS06}, respectively . 
Both definitions of 
$\mathrm{WF}_s$ and $\mathrm{HWF}_s$
are independent of the choice of $\mu>0.$
The usual $\mathrm{WF}_s$ is conic with respect to $\xi$. 
On the other hand, 
$\mathrm{HWF}_s$ is conic with respect to $(x, \xi).$ 

\vspace*{12pt}

We set 
\begin{equation}\label{eqn:hamiltonian}
p(x,\xi)=
\frac{1}{2}\sum_{j,k=1}^n a_{jk}(x)
\xi_j \xi_k
\end{equation}
and denote its Hamilton vector field by 
\[
H_p=\sum_{j=1}^n
\left(
\frac{\partial p}{\partial \xi_j}
\frac{\partial}{\partial x_j}
-\frac{\partial p}{\partial x_j}
\frac{\partial}{\partial \xi_j}
\right).
\]
Let $\gamma=\{(y(t),\eta(t)): t\in \mathbb R\}$
be an integral curve of $H_p$, that is, 
a solution to 
\begin{equation}
\dot{y}(t)
=\dfrac{\partial p}{\partial \xi}(y(t),\eta(t)), \quad
\dot{\eta}(t)
=-\dfrac{\partial p}{\partial x}(y(t),\eta(t)).
\end{equation}
We say that $\gamma$ is {\it backward nontrapping} 
if 
\begin{equation}\label{eqn:bnt}
 \lim_{t \to -\infty}|y(t)|=+\infty
\end{equation} 
holds. 
It is remarked that 
the nontrapping condition is necessary for 
some kind of smoothing effect. See Doi \cite{Do96} for the detail. 
We also remark that if $\gamma$ is 
backward nontrapping, then there exists the asymptotic
momentum
\begin{equation}\label{eqn:am}
\eta_-=\lim_{t \to -\infty} \eta(t) \in \mathbb R^n
\setminus \{0\} 
\end{equation}
under the assumption (A). 

We now state our main result. 
It is a analogy of 
\cite[Theorem 2.1]{MNS06} 
for the Gevrey case $s>1.$ 
\begin{theorem}\label{thm:phwf}
Assume that \textup{(A)} holds and that  
$\gamma$ is backward nontrapping. 
Let $\eta_-$ be the asymptotic momentum as 
$t \to -\infty.$ 
Assume that there exists $t_0>0$ such that 
\begin{equation}\label{eqn:phwf-1}
(-t_0 \eta_-,\eta_-) \notin \textup{HWF}_s(u_0), 
\end{equation} 
then we have 
\begin{equation}\label{eqn:phwf-2}
((t-t_0) \eta_-,\eta_-) \notin 
\textup{HWF}_s(u(t, \cdot))
\quad (0<t<\min (t_0, T)).
\end{equation}
Moreover, if $t_0<T$, then 
\begin{equation}\label{eqn:phwf-3}
\gamma \cap \textup{WF}_s (u(t,\cdot))=\emptyset
\end{equation}
holds for all $t$ close enough to $t_0$.
\end{theorem}
Theorem \ref{thm:phwf} means that
the microlocal Gevrey singularity of order $s$
appears only when the $\mathrm{HWF}_s$ hits $\{x=0\}$ 
as in the $C^{\infty}$ case or the analytic case. 
We remark that $(0, \eta_-) \notin \mathrm{HWF}_{s}(u)$ 
implies 
$(x, \eta_-) \notin \mathrm{WF}_{s}(u)$ 
for any 
$x \in \mathbb R^n.$

From this theorem, we obtain two results on microlocal 
Gevrey smoothing effects. 
To $\gamma=\{(y(t),\eta(t)): t\in \mathbb R\}$ 
and $\varepsilon>0$, we associate the set 
\[
 \Gamma_{\varepsilon}
=\bigcup_{t \le 0} \{
x\in \mathbb R^n : |x-y(t)|\le \varepsilon (1+|t|)
\}.
\]
The first one is concerned with the rapidly decaying data. 
\begin{corollary}\label{cor:decay}
Assume that \textup{(A)} holds and that
$\gamma$ is backward nontrapping. 
Assume that 
$e^{\delta_0 |x|^{1/s}}u_0 \in 
L^2 (\Gamma_{\varepsilon_0})$ 
for some $\delta_0>0$ and $\varepsilon_0>0$, 
then we have 
\begin{equation*}
\gamma \cap \textup{WF}_s (u(t,\cdot))=\emptyset
\quad (0<t<T).
\end{equation*}
\end{corollary}
We remark that in this case the condition 
\eqref{eqn:phwf-1} is satisfied for any 
$t_0>0.$

The second one deals with the initial data 
satisfying the mixed momentum condition, 
which also asserts the microlocal smallness of the data. 
It is a analogy of the result by Morimoto-Robbiano-Zuily \cite{MRZ}.
\begin{corollary}\label{cor:mix}
Assume that \textup{(A)} holds and that
$\gamma$ is backward nontrapping. 
Let $\eta_-$ be 
the asymptotic momentum as $t \to -\infty$. 
Assume that there exist
$\psi(\xi)\in C^{\infty}(\mathbb R^n)$
which equals to $1$ 
in a conic neighborhood of $\eta_-$, 
and $\varepsilon_0>0$, $A_0>0$, $A_1>0$ such that
\begin{equation*}
\|(x \cdot D_x)^l \psi(D_x)u_0\|
_{L^2(\Gamma_{\varepsilon_0})} 
\le A_0 A_1^l l!^{2s} \quad (l \in \mathbb N),
\end{equation*}
then we have 
\begin{equation*}
\gamma \cap \textup{WF}_s (u(t,\cdot))=\emptyset
\quad (0<t<T).
\end{equation*}
\end{corollary}
The paper is organized as follows. 
In Section \ref{sec:aae}
we recall the almost analytic extension of Gevrey symbols 
by Jung.
In Section \ref{sec:mee} 
we prove the microlocal exponential estimates 
in a Gevrey class. 
In Section \ref{sec:thm} and \ref{sec:cor}, 
we prove Theorem \ref{thm:phwf} 
and Corollary \ref{cor:mix}, respectively.

\vspace*{12pt}

\noindent{\bf Acknowledgment.} \enskip
The author would like to express his sincere gratitude 
to Hiroyuki Chihara for several discussions 
and valuable suggestions.

\section{Almost analytic extension}\label{sec:aae}

In this section we define the almost analytic extension 
of Gevrey functions following Jung's idea \cite{Ju}. 
We extend the functions on $\mathbb R^n$ to the complex strip 
with the parameterized width. 
We construct the extension to minimize 
the antiholomorphic derivatives of it 
on the complex strip.

Let $f(x)\in C^{\infty}(\mathbb R^n)$ 
satisfy the following condition.

\vspace*{12pt}

\noindent{\bf Assumption (G)}\enskip
There exist $s>1$, $C>0$, $R>0$ and $a \in \mathbb R$ such that
\begin{equation}
|\partial_x^{\alpha}f(x)|
\le CR^{|\alpha|}\alpha!^s \langle x \rangle^{a-|\alpha|}
\quad (\alpha \in \mathbb Z_+^n, x \in \mathbb R^n). 
\end{equation}

\vspace*{12pt}

\noindent
For $w>0$, we set the complex domain 
\[
S_w=\{z=x+iy \in \mathbb C^n : x \in \mathbb R^n, \quad|y_j|<w \quad
(1\le j \le n)\}. 
\]

\noindent
We define the almost analytic extension of $f$ on $S_w$ by 
\begin{equation}\label{def:aae}
\widetilde{f}(z)=\widetilde{f}(x+iy)
=\sum_{|\alpha|\le N(s,R,w)}\frac{(iy)^{\alpha}}{\alpha !}
\partial_x^{\alpha}f(x),  
\end{equation}
where $N(s,R,w)=[(Rw)^{-1/(s-1)}]$. 
We denote by $[c]$ the greatest integer not greater than $c \in \mathbb R$.
It is easy to see that $N \to +\infty$ as $s \to 1$ or 
$w \to 0$.
Although the definition \eqref{def:aae} depends on 
three parameters $s>1$ (Gevrey index), $R>0$ and 
$w>0$ (width of the strip), 
our interest is the case 
where 
$s$ and $R$ are fixed and $w$ tends to zero.
Indeed, we choose $w=\mathscr{O}(h^{1-1/s})\,(0<h \ll 1)$ 
in the argument of section \ref{sec:mee}. 

We see that if the function satisfies (G) then its derivatives
also satisfy the conditions similar to (G)
with the same constant $s>1$ and other 
$C'>0$, $R'>0$ and $a' \in \mathbb R.$
In other words, the index $s$ of (G) is stable under 
the differentiation. 

The first lemma is elementary but useful in our arguments.
It claims that not only $s$ but also $R$ is stable
in a slightly modified sense.
\begin{lemma}\label{lem:refine}
Let $f(x)$ be a function satisfying \textup{(G)}. 
Then for any $\beta \in \mathbb Z_+^n$, 
there exists $C>0$ such that 
\begin{equation}\label{eqn:refine-1}
|\partial_x^{\alpha+\beta}f(x)|
\le C(1+|\alpha|)^C R^{|\alpha|}
\alpha!^s \langle x \rangle^{a-|\alpha+\beta|}
\quad (\alpha \in \mathbb Z_+^n, x \in \mathbb R^n).  
\end{equation} 
\end{lemma}
\proof
It suffices to show that for any $\beta \in \mathbb Z_+^n$, 
there exists $C>0$ such that 
\begin{equation}\label{eqn:refine-2}
(\alpha+\beta)! \le C (1+|\alpha|)^C \alpha! \quad (\alpha \in 
\mathbb Z_+^n). 
\end{equation}
We use an induction over $|\beta|.$

In the case $\beta=e_j$, it is easy to see that
\begin{align*}
(\alpha+\beta)!
&=\alpha!(\alpha_j+1) \\
&\le (1+|\alpha|)\alpha!. 
\end{align*} 

We assume that \eqref{eqn:refine-2} holds up to $|\beta|\le m-1.$
In the case $|\beta|=m$, since $\beta_j \ge 1$ for some $j$,  
we have 
\begin{align*}
(\alpha+\beta)!&=(\alpha+e_j+\beta-e_j)! \\
&\le C(1+|\alpha+e_j|)^C (\alpha+e_j)! \\
&=C (1+|\alpha|)^C \left(\frac{2+|\alpha|}{1+|\alpha|}\right)^{C}
\alpha!(\alpha_j+1) \\
&\le 2^C C (1+|\alpha|)^{C+1} \alpha!,  
\end{align*}
which completes the proof.
\qed

\vspace*{12pt}

The behavior of the almost analytic extension 
as $w$ tends to $0$ 
is stated as follows. 
Especially, the proof of \eqref{eqn:aae-2}
implies that $N=N(s,R,w)$ is chosen to minimize 
$\overline{\partial}\widetilde{f}$ on $S_w.$
\begin{proposition}\label{prop:aae}
Let $f(x)$ be a function satisfying \textup{(G)} and 
$\widetilde{f}(x+iy)$ be its almost analytic extension 
on $S_w$ defined by \eqref{def:aae}. 
\begin{enumerate}
\item[\textup{(1)}]
For any $\alpha$, $\beta\in \mathbb Z_+^n$, 
there exists $C>0$
such that 
\begin{equation}\label{eqn:aae-1}
\sup_{x+iy \in S_w}
\frac{|
\partial_x^{\alpha}\partial_y^{\beta} \widetilde{f}(x+iy)|}
{\langle x \rangle^{a-|\alpha+\beta|}}
\le C 
\end{equation}
for $w \in (0, 1]$.
\item[\textup{(2)}]
For any $\alpha$, $\beta\in \mathbb Z_+^n$, 
there exist $C>0$ and $l>0$ such that 
\begin{equation}\label{eqn:aae-2}
\sup_{x+iy \in S_w}
\frac{
|\partial_x^{\alpha}\partial_y^{\beta} 
\overline{\partial}_j
\widetilde{f}(x+iy)|}
{\langle x \rangle^{a-|\alpha+\beta|-1}}
\le Cw^{-l}\exp\left(-\frac{\Omega}{w^{1/(s-1)}}\right) 
\end{equation}
for $w \in (0, 1]$, $1 \le j \le n$, 
where 
$\displaystyle\overline{\partial}_j=\frac{1}{2}
(\partial_{x_j}+i\partial_{y_j})$
and 
$\Omega=\dfrac{s-1}{R^{1/(s-1)}}.$
\end{enumerate} 
\end{proposition}

\proof

(1)\enskip We write $N$ instead of $N(s,R,w)$. 
By the definition \eqref{def:aae}, we have
\begin{align*}
\partial_x^{\alpha} \partial_y^{\beta} \widetilde{f}(x+iy)
&=\sum_{\substack{|\gamma| \le N  \\ \gamma \ge \beta}}
i^{|\gamma|} \frac {y^{\gamma-\beta}}{(\gamma-\beta)!}
\partial_x^{\gamma+\alpha} f(x) \\
&=\sum_{|\gamma| \le N-|\beta|}
i^{|\gamma+\beta|} \frac{y^{\gamma}}{\gamma!}
\partial_x^{\gamma+\alpha+\beta} f(x).  
\end{align*}  
Using $|y^{\gamma}|\le w^{|\gamma|}$, 
Lemma \ref{lem:refine} and 
\begin{equation}\label{eqn:stirling}
\alpha! \le c 
(1+|\alpha|)^{1/2}
\frac{|\alpha|^{|\alpha|}}{e^{|\alpha|}}
\quad (\alpha \in \mathbb Z_+^n) \quad \text{for some} \quad c>0, 
\end{equation}
which follows from Stirling's formula, we deduce that 
\begin{equation}\label{eqn:aae-3}
|\partial_x^{\alpha} \partial_y^{\beta} \widetilde{f}(x+iy)|
\le 
c_1 \langle x \rangle^{a-|\alpha+\beta|}
\sum_{|\gamma|\le N-|\beta|}
(1+|\gamma|)^{c_1} (Rw)^{|\gamma|}
\left(\frac{|\gamma|}{e}\right)^{|\gamma|(s-1)}.
\end{equation}
It remains to show that the sum in the
right hand side of \eqref{eqn:aae-3}
is uniformly bounded with respect to $w$. 
Set $M=N-|\beta|$. 
Then we have 
\begin{align*}
\text{(The sum in RHS of \eqref{eqn:aae-3})}
&=
\sum_{l=0}^M
(1+l)^{c_1} (Rw)^l
\left(\frac{l}{e}\right)^{l(s-1)}
\sum_{|\gamma|=l} 1 \\
&
\le  
\sum_{l=0}^M
(1+l)^{c_1 +n} e^{l(1-s)}, 
\end{align*}
where we use $\sum_{|\gamma|=l}1 \le (1+l)^n$ and 
$l\le M \le N \le (Rw)^{-1/(s-1)}.$

Pick up $L>1$ and $l_0 \in \mathbb N$ satisfying 
\[
 L e^{1-s}<1 \quad  \text{and} 
\quad
(1+l)^{c_1+n} \le L^l \quad (l > l_0).
\]
Then if $M \gg 1$, we obtain
\begin{align*}
\sum_{l=0}^M (1+l)^{c_1+n} e^{l(1-s)}
&\le 
\sum_{l=0}^{l_0} (1+l)^{c_1+n} e^{l(1-s)}
+\sum_{l=l_0 +1}^M (Le^{1-s})^l \\
&\le 
\sum_{l=0}^{l_0} (1+l)^{c_1+n} e^{l(1-s)}
+\frac{1}{1-Le^{1-s}}, 
\end{align*}
which implies the claim with $0< w \ll 1$. 
It is easy to see that the claim holds if $w$ is away from $0$.

(2)\enskip 
We first observe that 
\begin{equation*}
(\partial_{x_j}+i\partial_{y_j})\widetilde{f}(x+iy)
=\sum_{|\gamma|=N} \frac{(iy)^{\gamma}}{\gamma!}
\partial_{x}^{\gamma+e_j}f(x).
\end{equation*}
Then, just as in showing \eqref{eqn:aae-3}, 
we deduce that for any 
$\alpha, \beta \in \mathbb Z_+^n$, there exists $c_2>0$ such that
\begin{align*}
&|\partial_x^{\alpha}
\partial_y^{\beta}
(\partial_{x_j}+i\partial_{y_j})\widetilde{f}(x+iy)|  \\
&\le c_2 \langle x \rangle^{a-|\alpha+\beta|-1}
\sum_{|\gamma|= N-|\beta|}
(1+|\gamma|)^{c_2} (Rw)^{|\gamma|}
\left(\frac{|\gamma|}{e}\right)^{|\gamma|(s-1)}
\end{align*}
for $w\in (0,1]$ and $x+iy \in S_w$.
Set $M=N-|\beta|$. Then the sum in the right hand side
equals to 
\begin{align}\label{eqn:aae-4}
(1+M)^{c_2} \left(\frac{Rw M^{s-1}}{e^{s-1}}\right)^M
\sum_{|\gamma|=M} 1
\le (1+M)^{c_2 +n} \left(\frac{Rw M^{s-1}}{e^{s-1}}\right)^M.
\end{align}
A simple calculation shows that 
the function 
$a(x)=(Rwx^{s-1}/e^{s-1})^x\,(x>0)$ has the minimum
\[
 \exp \left(-\frac{\Omega}{w^{1/(s-1)}} \right)
\quad \left(\Omega=\frac{s-1}{R^{1/(s-1)}}\right)
\] 
at $x=(Rw)^{-1/(s-1)}$.
Since
$N=[(Rw)^{-1/(s-1)}]$  
and 
$M=N-|\beta| \sim (Rw)^{-1/(s-1)}$ as $w$ tends to $0$, 
we deduce that in the case of $0< w \ll 1$ there exist $c_3>0$ 
and $c_4>0$ such that
\[
 a(M) \le c_3 w^{-c_4} 
\exp \left(-\frac{\Omega}{w^{1/(s-1)}}\right).
\]
Combining this and \eqref{eqn:aae-4}, we obtain 
\eqref{eqn:aae-2} if $0< w \ll 1$.
We omit proof of the case $w$ away from zero.
\qed

\vspace*{12pt}

\noindent{\bf Remark.}\enskip 
We see that the right hand side of \eqref{eqn:aae-2} 
tends to $0$ as $s \to 1$ or $w \to 0.$
It is also clear that we can replace $w \in (0,1]$ in the statement 
of Proposition \ref{prop:aae} with $w \in (0,w_0]$ 
for arbitrary fixed $w_0>0.$

\section{Microlocal exponential estimates}\label{sec:mee}

This section is devoted to 
the microlocal exponential estimates
in a Gevrey class. 
Following \cite{MNS06}, we consider the estimates
with two parameters $h>0$ and $\mu>0.$
Roughly speaking, 
the former is a scaling parameter with respect to $\xi$, 
and the latter is a one to $x.$
We note that the analyticity of symbols, 
which we cannot use in our problem, is a essential
assumption in these estimates \cite{Ma}. 
To overcome this difficulty, 
we use the almost analytic extension of symbols defined in 
the previous section.
Moreover we assume the extra condition on the weight function. 

We first introduce a weight function $\psi$  
and a cutoff function $f$.

\vspace*{12pt}
\noindent{\bf Assumption (W1)}\enskip
Let $\psi(x,\xi)\in C^{\infty}(\mathbb R^{2n})$ 
be an $(h,\mu)$-dependent function satisfying 
\begin{itemize}
\item There exists $C_1>1$ such that
\[
 \mathrm{supp}\,[\psi]\subset
\left\{
(x,\xi) \in \mathbb R^{2n}:
\frac{1}{C_1}\le |\xi| \le C_1,\,
\frac{1}{C_1 \mu} \le \langle x \rangle \le \frac{C_1}{\mu}
\right\}
\] 
for $h, \mu \in (0,1].$
\item For any $\alpha, \beta \in \mathbb Z_+^n$ 
there exist $C_{\alpha \beta}>0$ such that 
\[
 |\partial_x^{\alpha}\partial_y^{\beta}\psi(x,\xi)|
\le C_{\alpha \beta}\mu^{|\alpha|}
\quad 
(x,\xi \in \mathbb R^n, h,\mu \in (0,1]).
\]
\item There exists $\nu>0$ such that
\begin{align*}
&\sup_{(x, \xi)\in \mathbb R^{2n}}|\partial_x \psi(x,\xi)|<\nu,  
\sup_{(x, \xi)\in \mathbb R^{2n}}|\partial_{\xi} \psi(x,\xi)|<\nu, 
\\
& \sup_{(x, \xi)\in \mathbb R^{2n}} |\psi(x,\xi)|
< \frac{s-1}{4(K_0 \nu)^{1/(s-1)}} 
\quad (h,\mu \in (0,1]),
\end{align*}
where $K_0$ is the same constant as in (A).
\end{itemize}

\vspace*{12pt}

\noindent{\bf Remark.}\enskip
The last assumption on the size of $\psi$
seems to be strong. 
However this extra condition enables us to 
treat new error terms particular to Gevrey cases
as negligible ones.

\vspace*{12pt}

\noindent{\bf Assumption (W2)}\enskip
Let $f(x,\xi)\in C^{\infty}(\mathbb R^{2n})$ 
be an $(h,\mu)$-dependent function satisfying 
\begin{itemize}
\item There exists $C_2>C_1$ such that
\[
 \mathrm{supp}\, [f] \subset
\left\{
(x,\xi) \in \mathbb R^{2n}:
\frac{1}{C_2}\le |\xi| \le C_2,\,
\frac{1}{C_2 \mu} \le \langle x \rangle \le \frac{C_2}{\mu}
\right\}
\] 
for $h, \mu \in (0,1].$
\item $f \equiv 1$ on $\mathrm{supp}\,[\psi].$
\item For any $\alpha, \beta \in \mathbb Z_+^n$ 
there exist $C_{\alpha \beta}>0$ such that 
\[
 |\partial_x^{\alpha}\partial_y^{\beta}f(x,\xi)|
\le C_{\alpha \beta}\mu^{|\alpha|}
\quad 
(x,\xi \in \mathbb R^n, h,\mu \in (0,1]).
\]
\item
$f \ge 0$ and $\sqrt{f} \in C^{\infty}(\mathbb R^{2n})$ 
satisfies the same estimates as above. 
\end{itemize}

\vspace*{12pt}

In the following we use the notation 
$\partial_{\mu}=\mu^{-1}\partial_x+i\partial_{\xi}.$ 
We set
\begin{align*}
&a(x,\xi)
=h^{-2} \frac{1}{2}\sum_{j,k=1}^n a_{jk}(x)\xi_j \xi_k
+h^{-1} \sum_{j=1}^n b_j(x)\xi_j +c(x), \\
&a_{\psi}(x,\xi)
=\widetilde{a}(x-h^{1-1/s}\partial_{\mu}\psi(x,\xi),
\xi+ih^{1-1/s}\mu \partial_{\mu}\psi(x,\xi)), 
\end{align*}
where $\widetilde{a}(x+iy, \xi)$ 
denotes the almost analytic extension
of $a(x, \xi)$ 
with respect to $x$
defined by \eqref{def:aae} with $w=h^{1-1/s}\nu$
and $R=K_0$, that is, 
\[
 \widetilde{a}(x+iy, \xi)
=\sum_{|\alpha| \le [ \nu_0 h^{-1/s}] }
\frac{(iy)^{\alpha}}{\alpha!}\partial_x^{\alpha}a(x, \xi)
\quad (x+iy \in S_{h^{1-1/s}\nu}, 
\xi \in \mathbb R^n)
\]
where
$\nu_0=(K_0 \nu)^{-1/(s-1)}.$
We note that
$a_{\psi}(x,\xi)$ is well-defined 
under the conditions (A) and (W1).   
Indeed, $a(x,\xi)$ is analytic in $\xi$ and 
\[
 |\mathrm{Im}\,(x-h^{1-1/s}\partial_{\mu}\psi(x,\xi))|
=h^{1-1/s} |\partial_{\xi}\psi(x,\xi)| <h^{1-1/s}\nu=w.
\]

Our main result in this section is the following. 
We write $T=T_{h,\mu}$ for simplicity.
\begin{theorem}\label{thm:mee}
Assume that \textup{(A)}, \textup{(W1)} and \textup{(W2)} hold. 
Suppose that there exists $d>0$ such that 
$0<h/\mu \le d$. 
Then there exists $C>0$ such that
\begin{align*}
&|\langle 
e^{\psi/h^{1/s}}Tu, fe^{\psi/h^{1/s}}TPu \rangle
-\langle e^{\psi/h^{1/s}}Tu, fa_{\psi}e^{\psi/h^{1/s}}Tu 
\rangle| \\
&\le 
C(h^{-1}\mu+\mu^{\sigma}+h \mu^{\sigma-1})
(\|\sqrt{f} e^{\psi/h^{1/s}}Tu \|^2+\|u\|^2) 
\end{align*}
for $u \in L^2(\mathbb R^n).$
\end{theorem}

To prove Theorem \ref{thm:mee}, we need some preliminaries.
We see that 
\begin{equation*}
P=h^{-2}\frac{1}{2}(hD_x)^2 
+h^{-2}p_2^W (x,hD_x)+h^{-1}p_1^W(x,hD_x)+p_0^W (x,hD_x),
\end{equation*}
where
\begin{align*}
&p_2(x,\xi)=\frac{1}{2}\sum_{j,k=1}^n (a_{jk}(x)-\delta_{jk})
\xi_j \xi_k, \\ 
&p_1(x,\xi)=\sum_{j=1}^n
b_j(x)\xi_j-\frac{1}{2i}
\sum_{j,k=1}^n \partial_{x_j}a_{jk}(x)\xi_k,\\
&p_0(x,\xi)
=p_0(x)
=c(x)
-\frac{1}{2i}
\sum_{j=1}^n
\partial_{x_j} b_j(x)
-\frac{1}{8}
\sum_{j,k=1}^n \partial_{x_j}\partial_{x_k}
a_{jk}(x)
\end{align*}
and $p_j^W(x,hD_x)$ denotes the Weyl-H\"ormander quantization
of $p_j$, that is 
\[
p_j^W(x,hD_x)u(x)
=\frac{1}{(2\pi h)^n}
\int \!\!\! \int_{\mathbb R^{2n}} 
e^{i(x-y)\cdot \xi/h}
p_j
\left(
\frac{x+y}{2}, \xi
\right)
u(y)\,dyd\xi
\]
for $u\in \mathscr{S}(\mathbb R^n).$ 
Hereafter, we use the $S(m,g)$ symbol class notation due to 
H\"ormander \cite{Ho}. We denote the set of operators
with their symbol in $S(m,g)$ by $\mathrm{OPS}(m,g).$ 

We set 
\[
q_j(x,\xi,x^{*},\xi^{*})=p_j(x-\xi^{*},x^{*})
\quad (x,\xi,x^{\ast},\xi^{*}\in \mathbb R^n)
\]
for $j=0,1,2.$ Then we have
\[
 TP_j=Q_j T \quad (j=0,1,2)
\]
where $P_j=p_j^W(x,hD_x)$ and 
$Q_j=q_j^W(x,\xi,hD_x,hD_{\xi})$. 
Let $g_1$ be a metric on $\mathbb R^{4n}$ defined by
\[
g_1=\frac {dx^2}{\Phi^2}
+\frac{d\xi^2}{\Psi^2}
+\frac{d{\xi^*}^2}{\Psi^2}
+\frac{d{\xi^*}^2}{\Phi^2}, 
\]
where
\begin{align*}
&\Phi(x,\xi,x^*,\xi^*)=\Phi(x,\xi^*)
=\sqrt{1+\frac{|x|^2}{\langle \xi^* \rangle^2}}, \\
&\Psi(x,\xi,x^*,\xi^*)=\Psi(\xi, x^*)
=\sqrt{1+\frac{|\xi|^2}{\langle \xi-x^* \rangle^2}}.   
\end{align*}
We recall that $q_j \in S(m_j, g_1)\,(j=0,1,2)$, where
\[
m_j
=
\langle \xi \rangle^j 
\langle x \rangle^{2-j-\sigma}
\langle \xi-x^* \rangle^j
\langle \xi^* \rangle^{|2-j-\sigma|}.
\]
See \cite[Lemma 3.6]{MNS06} for the details. 

We here give the fundamental lemma in the pseudodifferential 
calculus on the range of $e^{\psi/h^{1/s}}T$ 
without the proof. 
It follows from 
\cite[Lemma 3.5]{MNS06} 
and 
\cite[Lemma B.1]{MNS06}
with replacing $\psi \mapsto h^{1-1/s}\psi.$ 
\begin{lemma}\label{lem:form}
Assume \textup{(W1)}, \textup{(W2)}.
Suppose $Q \in \mathrm{OPS}
(\langle \xi \rangle^a
\langle x \rangle^b 
\langle \xi-x^* \rangle^m 
\langle \xi^* \rangle^l, g_1)$
with some $a,b,m,l \in \mathbb R.$ 
\begin{enumerate}
\item[\textup{(1)}] 
There exists $C>0$ such that 
\begin{align*}
|\langle e^{\psi/h^{1/s}}Tu, fQ e^{\psi/h^{1/s}}Tu \rangle|
\le 
C \mu^{-b}
(\|\sqrt{f}  e^{\psi/h^{1/s}} Tu\|^2+\|u\|^2 ) 
\end{align*} 
for $u \in \mathscr{S}(\mathbb R^n)$, 
$h,\mu \in (0,1]$ with $h/\mu \le d.$
\item[\textup{(2)}]
Suppose the symbol $Q$ has an asymptotic expansion 
supported in $\mathrm{supp}\,[\nabla f]$. 
Then for any $N>0$, 
there exists $C>0$ such that 
\begin{align*}
|\langle e^{\psi/h^{1/s}}Tu, Q e^{\psi/h^{1/s}}Tu \rangle|
\le 
C
(h^N \mu^N \|\sqrt{f}  e^{\psi/h^{1/s}} Tu\|^2+\mu^{-b} \|u\|^2 ) 
\end{align*} 
for $u \in \mathscr{S}(\mathbb R^n)$, 
$h,\mu \in (0,1]$ with $h/\mu \le d.$
\end{enumerate}
\end{lemma}
As mentioned in the above, 
under the assumption (A) 
we can define the almost analytic extension
of $p_j(x, \xi)$ with respect to $x$
with the setting $w=h^{1-1/s}\nu$ and $R=K_0$.  
We denote it by 
$\widetilde{p_j}(x+iy, \xi)\,
(x+iy \in S_{h^{1-1/s}\nu}, \xi \in \mathbb R^n)$. 
Therefore we can also define the almost analytic extension of 
$q_j(x,\xi,x^*,\xi^*)$ 
with respect to 
$\xi^*$ 
by
\begin{equation*}
 \widetilde{q_j}(x,\xi,x^*,\xi^*+i\eta^*)
=\widetilde{p_j}(x-\xi^*-i\eta^*, x^*)
\quad (x, \xi, x^* \in \mathbb R^n, 
\xi^*+i \eta^* \in S_{h^{1-1/s}\nu}).
\end{equation*}
We note that 
$\widetilde{q_j}$ 
is analytic with respect to
$x^*.$ 

Set
\[
 Q_{j \psi}=q_{j \psi}^W (x,\xi,hD_x,hD_{\xi})
\quad (j=0,1,2),
\]
where
\begin{align*}
q_{j \psi} (x, \xi, x^*, \xi^*)
&=\widetilde q_j (x, \xi, x^*+ih^{1-1/s}\partial_x \psi, 
\xi^*+ih^{1-1/s}\partial_{\xi} \psi) \\
&=\widetilde{p_j}
(x-\xi^*-ih^{1-1/s} \partial_{\xi} \psi, 
x^*+ih^{1-1/s}\partial_x \psi).
\end{align*}
Then we have the next key lemma, 
which should be compared to (3.5) in \cite{MNS06}, 
concerning the operator
\[ 
R_j=
e^{\psi/h^{1/s}} Q_j e^{-\psi/h^{1/s}}-Q_{j\psi}
\quad (j=0,1,2).
\]
It claims that, 
although the estimate of the remainder term becomes 
a slightly bad compared to the analytic case, 
we can approximate 
the distorted operator $e^{\psi/h^{1/s}} Q_j e^{-\psi/h^{1/s}}$
with $Q_{j\psi}$ in some sense. 
\begin{lemma}\label{lem:remainder}
Assume \textup{(W1)}, \textup{(W2)}. 
There exists $C>0$ such that
\begin{equation}\label{eqn:remainder-1}
|\langle 
e^{\psi/h^{1/s}}Tu, fR_j e^{\psi/h^{1/s}}Tu \rangle|
\le 
C h \mu^{\sigma+j-1}
(\|\sqrt{f}  e^{\psi/h^{1/s}} Tu\|^2+\|u\|^2 ) 
\end{equation} 
for $u \in \mathscr{S}(\mathbb R^n)$, $j=0,1,2$ 
and $h,\mu \in (0,1]$ with $h/\mu \le d.$ 
\end{lemma}
\proof
It suffices to show that
\begin{equation}\label{eqn:remainder-2}
\begin{split}
R_j \in \,
&\mathrm{OPS}\,(
h
\langle \xi \rangle^{j-1}
\langle x \rangle^{1-j-\sigma}
\langle \xi-x^* \rangle^{j+1}
\langle \xi^* \rangle^{|2-j-\sigma|+1}, 
g_1) \\
&+\mathrm{OPS}\,(
h^{\infty}
\langle \xi \rangle^{j}
\langle x \rangle^{1-j-\sigma}
\langle \xi-x^* \rangle^{j}
\langle \xi^* \rangle^{|1-j-\sigma|}, 
g_1).
\end{split}
\end{equation} 
Indeed, applying Lemma \ref{lem:form} (1) to $R_j$, 
we obtain \eqref{eqn:remainder-1}.

We denote by $X=(x,\xi)$, $Y=(y,\eta)$ and 
$X^*=(x^*, \xi^*)$ the points of $\mathbb R^{2n}.$ 
We see that 
\begin{align*}
&e^{\psi/h^{1/s}} Q_j e^{-\psi/h^{1/s}}u(X) \\
&=\frac{1}{(2\pi h)^{2n}}\int \!\!\! \int_{\mathbb R^{4n}}
e^{i(X-Y)\cdot X^*/h +(\psi(X)-\psi(Y))/h^{1/s}}
q_j \left(\frac{X+Y}{2}, X^*\right)u(Y)\,dYdX^* \\
&
=\frac{1}{(2\pi h)^{2n}} \int \!\!\! \int_{\mathbb R^{4n}} 
e^{i(X-Y) \cdot (X^*-ih^{1-1/s}\phi(X,Y))/h} 
q_j \left(\frac{X+Y}{2}, X^* \right) u(Y)\,dYdX^*  
\end{align*}
where
\begin{equation}\label{eqn:remainder-3}
\begin{split}
\phi(X,Y)&=(\phi_1(X,Y),\dots, \phi_{2n}(X,Y)), \\
\phi_l(X,Y)&=
\int_0^1 \partial_{X_l} 
\psi(Y_1,\dots, Y_{l-1}, Y_l +\tau (X_l-Y_l), X_{l+1}, \dots, 
X_{2n})\,d\tau.
\end{split}
\end{equation}
We remark that our choice of $\phi_l$ is not a standard one.
For instance, $\phi_l (X,Y) \ne \phi_l(Y,X)$. 
However our $\phi_l$  
satisfies
\begin{equation}\label{eqn:remainder-4}
\begin{split}
&(X_l-Y_l)\phi_l(X,Y) \\
&=\psi(Y_1,\dots, Y_{l-1}, X_l, \dots, X_{2n})- 
\psi(Y_1,\dots, Y_{l}, X_{l+1}, \dots, X_{2n}). 
\end{split}
\end{equation}
In particular, $(X_l-Y_l)\phi_l(X,Y)$ is bounded
on $\mathbb R^{4n}$ 
which is useful to estimate some error terms 
particular to the Gevrey case.

We here change the integral variables 
$X_1^* \mapsto X_1^*+ih^{1-1/s}\phi_1$ 
up to 
$X_{2n}^* \mapsto X_{2n}^*+ih^{1-1/s}\phi_{2n}$ 
successively
in the above expression of 
$e^{\psi/h^{1/s}} Q_j e^{-\psi/h^{1/s}}u$.
Applying the Stokes formula, we have
\begin{align*}
&e^{\psi/h^{1/s}} Q_j e^{-\psi/h^{1/s}}u(X) \\
&
=\frac{1}{(2\pi h)^{2n}} 
\int \!\!\! \int_{\mathbb R^{4n}} 
e^{i(X-Y) \cdot X^* /h} 
\widetilde{q_j} \left( \frac{X+Y}{2}, X^*+ih^{1-1/s} \phi(X,Y) \right) 
u(Y) \, dYdX^* \\
&\qquad{}+\sum_{k=n+1}^{2n}
\frac{1}{(2\pi h)^{2n}} 
\int \!\!\! \int_{\mathbb R^{4n}} 
e^{i(X-Y) \cdot X^* /h}
r_{j,k}(X,Y,X^*) 
u(Y) \, dYdX^* \\
&=Q_j^{(1)}u+Q_j^{(2)}u, 
\end{align*}
where 
\begin{equation}\label{eqn:remainder-5}
\begin{split}
&r_{j,k}(X,Y,X^*) \\
&=2ih^{1-1/s}\phi_k (X, Y)
\int_0^1 e^{(X-Y) \cdot 
(0,\dots, 0, (1-\theta)\phi_k, \phi_{k+1}, \dots, \phi_{2n})/h^{1/s}}
\\
&\quad {} \times
\frac{\partial \widetilde{q_j}}{\partial \bar{Z_k^*}}
\left(\frac{X+Y}{2}, X^* +ih^{1-1/s}
(\phi_1, \dots, \phi_{k-1}, \theta \phi_k, 0, \dots, 0)
\right)\,d\theta
\end{split} 
\end{equation}
and
$Z^*=(z^*, \zeta^*)=X^* +iY^* =(x^*+iy^*, \xi^*+i\eta^*)
\in \mathbb C^{2n}.$
It is remarked that $r_{j,k}\equiv 0$ for $1 \le k \le n$, 
since
$\widetilde{q_j}$ is analytic with respect to 
$Z_k^* \, (1 \le k \le n).$
In the following we shall simplify the double symbols
$\widetilde{q_j} ((X+Y)/2, X^*+ih^{1-1/s} \phi(X,Y))$
and 
$r_{j,k}(X,Y,X^*).$

First we observe that for any 
$\alpha, \beta, \gamma, \delta \in \mathbb Z_+^n$, 
there exists $C>0$ such that 
\begin{equation}\label{eqn:remainder-6}
|\partial_x^{\alpha}\partial_{\xi}^{\beta}
\partial_y^{\gamma}\partial_{\eta}^{\delta}
\phi_l (x,\xi,y,\eta)|
\le 
C \frac{
\langle x-y \rangle^{|\alpha+\gamma|}
\langle \xi-\eta \rangle^{|\beta+\delta|}
}
{
\langle x \rangle^{|\alpha|}
\langle \xi \rangle^{|\beta|}
\langle y \rangle^{|\gamma|}
\langle \eta \rangle^{|\delta|}
} 
\end{equation}
$(x,\xi,y,\eta \in \mathbb R^n, 1 \le l \le 2n).$
Indeed, combining \eqref{eqn:remainder-3}, 
\[
|\partial_x^{\alpha}\partial_{\xi}^{\beta}\psi(x,\xi)|
\le 
C \langle x \rangle^{-|\alpha|}
\langle \xi \rangle^{-|\beta|}
\]
$(x,\xi \in 
\mathbb R^n, h,\mu \in (0, 1]))$, 
which follows from (W1), and 
\[
 \langle x \rangle \langle y \rangle
\le
C \langle x-y \rangle
\langle 
(y_1, \dots, y_{j-1}, y_j +\tau (x_j -y_j), 
x_{j+1}, \dots, x_n
) \rangle
\]
$(x, y \in \mathbb R^n, \tau \in [0, 1], 1 \le j \le n)$, 
we deduce \eqref{eqn:remainder-6} for 
$1 \le l \le n.$
The case 
$n+1 \le l \le 2n$ 
can be handled in the same manner. 
We note that $\phi_l\,(1 \le l \le n)$ is independent of $\eta$, 
and 
$\phi_l\,(n+1 \le l \le 2n)$ is independent of $x$. 

Applying Proposition \ref{prop:aae} (1) to 
$\widetilde{q_j}$, 
we see that for any 
$\alpha, \beta, \alpha^*, \beta^*, 
\gamma^*, \delta^* \in \mathbb Z_+^n$, 
there exists $C>0$ such that 
\begin{align*}
&|
\partial_x^{\alpha} \partial_{\xi}^{\beta}
\partial_{x^*}^{\alpha^*} \partial_{\xi^*}^{\beta^*}
\partial_{y^*}^{\gamma^*} \partial_{\eta^*}^{\delta^*}
\widetilde{q_j} (x,\xi, x^*+iy^*, \xi^*+i\eta^*)| \\
&\le C
\langle x-\xi^* \rangle^{2-j-\sigma-|\alpha+\beta^*+\delta^*|}
\langle x^* \rangle^{j-|\beta+\alpha^*+\gamma^*|} \\
&\le
C \frac
{
\langle \xi \rangle^{j}
\langle x \rangle^{2-j-\sigma}
\langle \xi-x^* \rangle^{j}
\langle \xi^* \rangle^{|2-j-\sigma|}
}
{\Phi(x,\xi^*)^{|\alpha+\beta^*+\delta^*|}
\Psi(\xi, x^*)^{|\beta+\alpha^*+\gamma^*|}}
\end{align*}
$(x,\xi \in \mathbb R^n, 
x^*+iy^*, \xi^*+i\eta^* \in S_{h^{1-1/s} \nu}, h \in (0, 1]
).$
We write 
\begin{equation*}
\phi(X, Y)=(\phi^{(1)}(X, Y), \phi^{(2)}(X, Y)),
\end{equation*}
that is 
\begin{align*}
&\phi^{(1)}(x,\xi,y,\eta)
=\phi^{(1)}(x,\xi,y)
=(\phi_1(x,\xi,y), \dots, \phi_n(x,\xi,y)), \\
&\phi^{(2)}(x,\xi,y,\eta)
=\phi^{(2)}(\xi,y,\eta)
=(\phi_{n+1}(\xi,y,\eta), \dots, \phi_{2n}(\xi,y,\eta)).
\end{align*}
Combining the above estimate of 
$\widetilde{q_j}$
and \eqref{eqn:remainder-6}, we deduce that
for any 
$\alpha, \beta, \gamma, \delta, \alpha^*, \beta^* 
\in \mathbb Z_+^n$, 
there exists 
$C>0$ 
such that 
\begin{equation}\label{eqn:remainder-7}
\begin{split}
&
\left|
\partial_x^{\alpha} \partial_{\xi}^{\beta}
\partial_y^{\gamma} \partial_{\eta}^{\delta}
\partial_{x^*}^{\alpha^*} \partial_{\xi^*}^{\beta^*}
\widetilde{q_j} 
\left(
\frac{x+y}{2},\frac{\xi+\eta}{2}, 
x^*+ih^{1-1/s}\phi^{(1)}, \xi^*+ih^{1-1/s}\phi^{(2)}
\right)
\right|
\\
&\le 
C 
\langle \xi+\eta \rangle^{j}
\langle x+y \rangle^{2-j-\sigma}
\langle x^*-(\xi+\eta)/2 \rangle^{j}
\langle \xi^* \rangle^{|2-j-\sigma|}
\\
&\quad {} \times 
\frac{
\langle x-y \rangle^{|\alpha+\gamma|}
\langle \xi-\eta \rangle^{|\beta+\delta|}
}
{\Phi(x,\xi^*)^{|\alpha|}
\Phi(y,\xi^*)^{|\gamma|}
\Phi(\frac{x+y}{2},\xi^*)^{|\beta^*|}
\Psi(\xi, x^*)^{|\beta|}
\Psi(\eta, x^*)^{|\delta|}
\Psi(\frac{\xi+\eta}{2}, x^*)^{|\alpha^*|}}
\end{split}
\end{equation}
$(x,\xi,y,\eta,x^*,\xi^* \in \mathbb R^n, h,\mu \in (0,1]).$
It follows from \eqref{eqn:remainder-7} 
that there exists a simplified symbol
\begin{equation}
\label{eqn:remainder-8}
\rho_j (x,\xi,x^*,\xi^*)\in S(m_j,g_1) 
\end{equation}
such that 
\begin{gather}\label{eqn:remainder-9}
\rho_j^W(x,\xi,hD_x,hD_{\xi})=Q_j^{(1)}, \\
\label{eqn:remainder-10}
\rho_j (x,\xi,x^*,\xi^*)
-q_{j \psi}(x,\xi,x^*,\xi^*)
\in S(h \Phi^{-1} \Psi^{-1} m_j, g_1).
\end{gather}
It is remarked that we cannot replace $h \Phi^{-1} \Psi^{-1}$
in \eqref{eqn:remainder-10} by $h^2 \Phi^{-2} \Psi^{-2}$
because of the property $\phi(X,Y) \ne \phi(Y,X).$

We next give the estimate of 
$r_{j,k}(X,Y,X^*)\,(n+1 \le k \le 2n).$
Using \eqref{eqn:remainder-4}, we rewrite \eqref{eqn:remainder-5}
as 
\begin{equation}\label{eqn:remainder-11}
\begin{split}
&r_{j,n+l}(x,\xi,y,\eta,x^*,\xi^*) \\
&=2ih^{1-1/s}\phi_{n+l}(x,\xi, y, \eta, x^*, \xi^*) \\
&\quad \times 
\int_0^1 
e^{
\{
(1-\theta)\psi(y,\eta_1, \dots, \eta_{l-1}, \xi_l, \dots, \xi_n)
+\theta \psi (y, \eta_1, \dots, \eta_l, \xi_{l+1}, \dots, \xi_n)
-\psi(y, \eta)\}/h^{1/s}
} \\
&\quad
\times
\frac{\partial \widetilde{q_j}}{\partial \bar{\zeta_l^*}}
\left(
\frac{x+y}{2}, \frac{\xi+\eta}{2}, 
x^*+ih^{1-\frac{1}{s}}\phi^{(1)}, 
\xi^*+ih^{1-\frac{1}{s}}
(\phi_{n+1}, \dots, \phi_{n+l-1}, \theta \phi_{n+l},
0,\dots, 0) 
\right) \\
&
\quad
\times d\theta \quad (1 \le l \le n).
\end{split}
\end{equation}
Applying Proposition \ref{prop:aae} (2), 
we see that for any 
$\alpha, \beta, \alpha^*, \beta^*, \gamma^*, \delta^* 
\in \mathbb Z_+^n$, 
there exists 
$C>0$ such that 
\begin{equation}\label{eqn:remainder-12}
\begin{split}
&
\left|
\partial_x^{\alpha} \partial_{\xi}^{\beta}
\partial_{x^*}^{\alpha^*} \partial_{\xi^*}^{\beta^*}
\partial_{y^*}^{\gamma^*} \partial_{\eta^*}^{\delta^*}
\overline{\partial}_{\zeta_l^*}
\widetilde{q_j}(x, \xi, x^*+iy^*, \xi^*+i\eta^*)
\right|
\\
&
\le 
C h^{-C} e^{-\Sigma_0/h^{1/s}}
\frac{
\langle \xi \rangle^{j}
\langle x \rangle^{1-j-\sigma}
\langle x^*-\xi \rangle^{j}
\langle \xi^* \rangle^{|1-j-\sigma|}
}
{\Phi(x,\xi^*)^{|\alpha+\beta^*+\delta^*|}
\Psi(\xi, x^*)^{|\beta+\alpha^*+\gamma^*|}}
\end{split}
\end{equation}
$(x,\xi \in \mathbb R^n,
x^*+iy^*, \xi^*+i\eta^* \in S_{h^{1-1/s}\nu}, 
h \in (0,1], 1 \le l \le n)$, 
where 
$\Sigma_0=\dfrac{s-1}{(K_0 \nu)^{1/(s-1)}}$
and 
$\overline{\partial}_{\zeta_l^*}
=(\partial_{\xi_l^*}+i\partial_{\eta_l^*})/2.$
We also see that 
\begin{equation}\label{eqn:remainder-13}
|\partial_x^{\alpha}
\partial_{\xi}^{\beta}
\partial_y^{\gamma}
\partial_{\eta}^{\delta}
\psi(Y_1, \dots, Y_{l-1}, X_l, \dots, X_{2n})|
\le C
\frac{
\langle x-y \rangle^{|\alpha+\gamma|}
\langle \xi-\eta \rangle^{|\beta+\delta|}
}
{
\langle x \rangle^{|\alpha|}
\langle \xi \rangle^{|\beta|}
\langle y \rangle^{|\gamma|}
\langle \eta \rangle^{|\delta|}}
\end{equation}
$(x,\xi,y,\eta \in \mathbb R^n, 
h, \mu \in (0,1], 
1 \le l \le 2n)$ 
and, using (W1),
\begin{equation}\label{eqn:remainder-14}
\begin{split}
\left| \,\bigl(
\text{exponential term in \eqref{eqn:remainder-11}}
\bigr) \,
\right|
&
\le \exp 
\left(
3 \sup_{(x, \xi) \in \mathbb R^{2n}}|\psi(x,\xi)| /
h^{1/s}
\right)
\\
&\le \exp \left(
\frac{3 \Sigma_0}{4h^{1/s}}
\right)
\end{split} 
\end{equation}
$(y,\eta,\xi \in \mathbb R^n, \theta \in [0,1], h,\mu \in (0,1]).$
Combining 
\eqref{eqn:remainder-6}, \eqref{eqn:remainder-12},
\eqref{eqn:remainder-13}, and \eqref{eqn:remainder-14}, 
we deduce that 
for any 
$\alpha, \beta, \gamma, \delta, \alpha^*, \beta^* 
\in \mathbb Z_+^n$, 
there exists 
$C>0$ such that 
\begin{equation}\label{eqn:remainder-15}
\begin{split}
&
\left|
\partial_x^{\alpha} \partial_{\xi}^{\beta}
\partial_{y}^{\gamma} \partial_{\eta}^{\delta}
\partial_{x^*}^{\alpha^*} \partial_{\xi^*}^{\beta^*}
r_{j,n+l} (x, \xi, y, \eta, x^*, \xi^*)
\right|
\\
&
\le 
C h^{-C} \exp \left(-\frac{\Sigma_0}{4h^{1/s}} 
\right)
\langle \xi+\eta \rangle^{j}
\langle x+y \rangle^{1-j-\sigma}
\langle x^*-(\xi+\eta)/2 \rangle^{j}
\langle \xi^* \rangle^{|1-j-\sigma|} \\
&\quad {} \times 
\frac{
\langle x-y \rangle^{|\alpha+\gamma|}
\langle \xi-\eta \rangle^{|\beta+\delta|}
}
{\Phi(x,\xi^*)^{|\alpha|}
\Phi(y,\xi^*)^{|\gamma|}
\Phi(\frac{x+y}{2},\xi^*)^{|\beta^*|}
\Psi(\xi, x^*)^{|\beta|}
\Psi(\eta, x^*)^{|\delta|}
\Psi(\frac{\xi+\eta}{2}, x^*)^{|\alpha^*|}}
\end{split}
\end{equation}
$(x,\xi, y, \eta,x^*, \xi^* \in \mathbb R^n,
h, \mu  \in (0,1], 1 \le l \le n)$.
It follows from \eqref{eqn:remainder-15} that 
there exists a simplified symbol
\begin{equation}\label{eqn:remainder-16}
r_j(x,\xi,x^*,\xi^*)\in 
S(h^{\infty}
\langle \xi \rangle^{j}
\langle x \rangle^{1-j-\sigma}
\langle x^*-\xi \rangle^{j}
\langle \xi^* \rangle^{|1-j-\sigma|}, g_1)
\end{equation}
such that 
\begin{equation}\label{eqn:remainder-17}
r_j^W (x, \xi, h D_x, h D_{\xi})=Q_j^{(2)}. 
\end{equation}
Then it follows from \eqref{eqn:remainder-8},
\eqref{eqn:remainder-9}, \eqref{eqn:remainder-10},
\eqref{eqn:remainder-16}
and 
\eqref{eqn:remainder-17}
that we obtain \eqref{eqn:remainder-2}, 
which completes the proof.
\qed

\vspace*{12pt}

\noindent{\bf Remark.}\enskip 
It is obvious that we can show a stronger result rather than 
\eqref{eqn:remainder-2}
in which $\mathscr{O}(h^{\infty})$ 
is replaced with 
$\mathscr{O}(e^{-c/h^{1/s}}).$

\vspace*{12pt}

To estimate 
$\langle 
e^{\psi/h^{1/s}}Tu, f Q_{j\psi}e^{\psi/h^{1/s}}Tu \rangle$, 
we use the following lemma. Set
\begin{align*}
p_{j\psi}(x,\xi)
&=q_{j \psi}(x,\xi,\xi-h^{1-1/s}\mu \partial_{\xi}\psi, 
h^{1-1/s}\mu^{-1} \partial_{x}\psi ) \\
&=\widetilde{p_j}
(x-h^{1-1/s} \partial_{\mu} \psi,
\xi+ih^{1-1/s} \mu \partial_{\mu} \psi)
\quad (j=0,1,2).
\end{align*}
\begin{lemma}\label{lem:conjugate}
Assume \textup{(W1)}, \textup{(W2)}. 
\begin{enumerate}
\item[\textup{(1)}] There exists $C>0$ such that 
\begin{align*}
&\left|\left\langle e^{\psi/h^{1/s}}Tu, 
f e^{\psi/h^{1/s}}
\left\{
\frac{1}{2}(hD_x)^2-\frac{1}{2}
( \xi+ih^{1-1/s}\mu \partial_{\mu}\psi(x,\xi))^2
+\frac{n}{4}h\mu \right\}
Tu 
\right\rangle \right| \\
&\le C h \mu 
(h^{1-1/s}\mu \|
\sqrt{f}e^{\psi/h^{1/s}}Tu\|^2
+\|u\|^2)
\end{align*} 
for $u \in \mathscr{S}(\mathbb R^n)$, 
$h, \mu \in (0, 1].$
\item[\textup{(2)}] There exists $C>0$ such that 
\begin{align*}
&|\langle e^{\psi/h^{1/s}}Tu, 
f \{Q_{j\psi}-p_{j \psi}(x, \xi) \}
e^{\psi/h^{1/s}} Tu 
\rangle | \\
&\le C h \mu^{\sigma+j-1} 
\| \sqrt{f}e^{\psi/h^{1/s}}Tu\|^2+\|u\|^2)
\end{align*} 
\end{enumerate}
for $u \in \mathscr{S}(\mathbb R^n)$, $j=0,1,2$
and $h, \mu \in (0, 1]$ with $h/\mu \le d.$
\end{lemma}
\proof 
(1)\enskip 
Combining the formula 
\begin{align*}
&\langle e^{g/h}Tu, f e^{g/h}(hD_x)^{\alpha}Tu \rangle \\
&
=\left\langle
 e^{g/h}Tu, 
\left\{\left(
\xi+i\mu \partial_{\mu}g+\frac{i}{2}h\mu \partial_{\mu}
\right)^{\alpha} f \right\} e^{g/h}(hD_x)^{\alpha}Tu
\right\rangle \quad (\alpha \in \mathbb Z_+^n)
\end{align*}
(\cite[Lemma C.1]{MNS06}) with the setting $g=h^{1-1/s}\psi$ and
\begin{align*}
&\left(\xi+ih^{1-1/s} \mu \partial_{\mu}\psi
+\frac{i}{2}h\mu \partial_{\mu} \right)^2 f \\
&=
(\xi+ih^{1-1/s}\mu \partial_{\mu}\psi)^2 f
-\frac{n}{2}h\mu f
-\frac{1}{2}h^{1-1/s}\mu^2 (\partial_{\mu}^2 \psi)f \\
&\quad{}
+ih\mu (\xi+ih^{1-1/s}\mu \partial_{\mu}\psi)\partial_{\mu}f
-\frac{1}{4}h^2 \mu^2 \partial_{\mu}^2 f,
\end{align*}
we deduce the claim.

(2)\enskip For $j=0,1,2$, we set 
\begin{align*}
&q_{j\psi}(x,\xi,x^*,\xi^*)
-q_{j\psi}(x, \xi, \xi-h^{1-1/s}\mu \partial_{\xi}\psi, 
h^{1-1/s}\mu^{-1} \partial_x \psi) \\ 
&=q_{j\psi}^{(1)}(x,\xi,x^*,\xi^*)
(x^*-\xi+h^{1-1/s}\mu \partial_{\xi}\psi)
+q_{j\psi}^{(2)}(x,\xi,x^*,\xi^*)
(\xi^*-h^{1-1/s}\mu^{-1} \partial_x \psi),
\end{align*}
where 
\begin{align*}
&q_{j\psi}^{(1)}(x,\xi,x^*,\xi^*) \\
&=
\int_0^1 \frac{\partial q_{j\psi}}{\partial x^*}
(x,
\xi, 
\theta x^*+(1-\theta)(\xi-h^{1-1/s}\mu \partial_{\xi}\psi),
\theta \xi^*+(1-\theta)h^{1-1/s}\mu^{-1} \partial_x \psi)
\,d\theta, \\
&q_{j\psi}^{(2)}(x,\xi,x^*,\xi^*) \\
&=
\int_0^1 \frac{\partial q_{j \psi}}{\partial \xi^*}
(x,
\xi, 
\theta x^*+(1-\theta)(\xi-h^{1-1/s}\mu \partial_{\xi}\psi),
\theta \xi^*+(1-\theta)h^{1-1/s}\mu^{-1} \partial_x \psi)
\,d\theta.
\end{align*}
Since $q_{j \psi}(x,\xi, x^*, \xi^*)\in S(m_j, g_1)$,
we can verify that 
\begin{align*}
&q_{j \psi}^{(1)} (x,\xi, x^*, \xi^*)
\in 
S(
\langle \xi \rangle^{j-1}
\langle x \rangle^{2-j-\sigma}
\langle \xi-x^* \rangle^{j+1}
\langle \xi^* \rangle^{|2-j-\sigma|}
, g_1), \\ 
&q_{j \psi}^{(2)} (x,\xi, x^*, \xi^*)
\in 
S(
\langle \xi \rangle^{j}
\langle x \rangle^{1-j-\sigma}
\langle \xi-x^* \rangle^{j}
\langle \xi^* \rangle^{|2-j-\sigma|+1}
, g_1).
\end{align*}
We also set
\begin{align*}
&A=h D_x-\xi+h^{1-1/s}\mu \partial_{\xi}\psi(x,\xi), \\
&B=h D_{\xi}-h^{1-1/s}\mu^{-1} \partial_{x} \psi(x,\xi)
\end{align*}
and 
\begin{align*}
Q_{j \psi}^{(k)}&=q_{j \psi}^{(k)W}(x,\xi,hD_x,hD_{\xi}) 
\quad (k=1,2),\\
R^{(j)}&=Q_{j\psi}
-q_{j\psi}
(x,\xi,\xi-h^{1-1/s}\mu \partial_{\xi}\psi, 
h^{1-1/s}\mu^{-1} \partial_{x}\psi) \\
&\qquad{}-\frac{1}{2}\{
A Q_{j \psi}^{(1)}+Q_{j \psi}^{(1)}A
+B Q_{j \psi}^{(2)}+Q_{j \psi}^{(2)}B
\}.
\end{align*}
Then, using the symbolic calculus, we deduce that 
\begin{equation}\label{eqn:conjugate-1}
R^{(j)}\in \mathrm{OPS}(h^2
\langle \xi \rangle^{j-2}
\langle x \rangle^{-j-\sigma}
\langle \xi-x^* \rangle^{j+3} 
\langle \xi^* \rangle^{|2-j-\sigma|+3}, g_1 
). 
\end{equation}
Especially, applying Lemma \ref{lem:form} (1), we obtain 
\begin{equation}\label{eqn:conjugate-2}
\begin{split}
&|\langle e^{\psi/h^{1/s}}Tu, 
f R^{(j)} e^{\psi/h^{1/s}} Tu \rangle | \\
&\le C h^2 \mu^{\sigma+j} 
\| \sqrt{f}e^{\psi/h^{1/s}}Tu\|^2+\|u\|^2)\quad (j=0,1,2).
\end{split}
\end{equation}

Next we shall show that 
\begin{equation}\label{eqn:conjugate-3}
\begin{split}
&|\langle e^{\psi/h^{1/s}}Tu, 
f (A Q_{j \psi}^{(1)}+Q_{j \psi}^{(1)}A
+B Q_{j \psi}^{(2)}+Q_{j \psi}^{(2)}B)
e^{\psi/h^{1/s}} Tu \rangle | \\
&\le C h \mu^{\sigma+j-1} 
\| \sqrt{f}e^{\psi/h^{1/s}}Tu\|^2+\|u\|^2)\quad (j=0,1,2).
\end{split}
\end{equation}
It follows from 
\eqref{eqn:conjugate-2} and \eqref{eqn:conjugate-3} that
we deduce the claim of (2). 

We write $T_{\psi}=e^{\psi/h^{1/s}}T$ for abbreviation. 
Since 
\[
 (A-i\mu B)T_{\psi}u=0,
\]
we see that 
\begin{align*}
T_{\psi}^* f (A Q_{j \psi}^{(1)}+Q_{j \psi}^{(1)}A
+B Q_{j \psi}^{(2)}+Q_{j \psi}^{(2)}B) T_{\psi}
=T_{\psi}^* (L_j^{(1)} +L_j^{(2)}) T_{\psi},
\end{align*}
where
\begin{align*}
&L_j^{(1)}
=if \bigl(-\mu \bigl[B, Q_{j\psi}^{(1)}\bigr]
+\mu^{-1}\bigl[ A, Q_{j\psi}^{(2)}\bigr]\bigr),\\
&L_j^{(2)}
=[f, A+i\mu B] Q_{j\psi}^{(1)}
-i\mu^{-1}[f, A+i\mu B] Q_{j\psi}^{(2)}.
\end{align*}
We also see that 
\begin{align*}
&\bigl[
B, Q_{j \psi}^{(1)}
\bigr]
\in \mathrm{OPS}(
h 
\langle \xi \rangle^{j-2}
\langle x \rangle^{2-j-\sigma}
\langle \xi-x^* \rangle^{j+2}
\langle \xi^* \rangle^{|2-j-\sigma|}
, g_1), \\
&\bigl[
A, Q_{j \psi}^{(2)}
\bigr] 
\in \mathrm{OPS}(
h 
\langle \xi \rangle^{j}
\langle x \rangle^{-j-\sigma}
\langle \xi-x^* \rangle^{j}
\langle \xi^* \rangle^{|2-j-\sigma|+2}
, g_1).
\end{align*}
Combining this and Lemma \ref{lem:form} (1), we have
\begin{equation}\label{eqn:conjugate-4}
\begin{split}
&|\langle e^{\psi/h^{1/s}} Tu, 
L_j^{(1)} e^{\psi/h^{1/s}}Tu
\rangle | \\
&\le C h \mu^{\sigma+j-1} 
(\|\sqrt{f}e^{\psi/h^{1/s}}Tu\|^2+\|u\|^2)
\quad (j=0,1,2).
\end{split} 
\end{equation}
It also follows from Lemma \ref{lem:form} (2) that we obtain
\begin{equation}\label{eqn:conjugate-5}
\begin{split}
&|\langle e^{\psi/h^{1/s}} Tu, 
L_j^{(2)} e^{\psi/h^{1/s}}Tu
\rangle | \\
&\le C h \mu^{\sigma+j-1} 
(\|\sqrt{f}e^{\psi/h^{1/s}}Tu\|^2+\|u\|^2)
\quad (j=0,1,2). 
\end{split} 
\end{equation}
Combining \eqref{eqn:conjugate-4} and \eqref{eqn:conjugate-5}, 
we deduce \eqref{eqn:conjugate-3}
and which completes the proof.
\qed
\vspace*{12pt}

We are now ready to prove Theorem \ref{thm:mee}. 

\vspace*{12pt}

\noindent{\it Proof of Theorem \textup{\ref{thm:mee}}.}\enskip
By the definition of $P_j$ and 
the relation
\[
TP_j=Q_j T
=e^{-\psi/h^{1/s}}
(Q_{j\psi}+R_j)
e^{\psi/h^{1/s}}T, 
\]
we have
\begin{align*}
&\langle e^{\psi/h^{1/s}}Tu, 
f e^{\psi/h^{1/s}} TPu \rangle   \\
&=
\frac{1}{2}h^{-2} \langle e^{\psi/h^{1/s}}Tu, 
f e^{\psi/h^{1/s}} (hD_x)^2 Tu \rangle   \\
&
\quad{}+\sum_{j=0}^2 h^{-j}
\langle e^{\psi/h^{1/s}}Tu, 
f (Q_{j \psi}+R_j) e^{\psi/h^{1/s}} Tu \rangle.  
\end{align*}
We note that 
\begin{align*}
&\left|
\frac{1}{2}h^{-2}(\xi+ih^{1-1/s}\mu\partial_{\mu}\psi(x,\xi))^2
+\sum_{j=0}^2 h^{-j} p_{j\psi}(x,\xi)-a_{\psi}(x,\xi)
\right| \\
& \le 
C (h^{-1}\langle x \rangle^{-\sigma-1}|\xi|+
\langle x \rangle^{-\sigma}).
\end{align*}
Combining this, Lemma \ref{lem:remainder} and 
Lemma \ref{lem:conjugate}, 
we obtain our claim.
\qed

\vspace*{12pt}

Finally, we give the useful corollary of Theorem \ref{thm:mee}.
\begin{lemma}\label{lem:taylor}
Under the same assumptions as in Theorem \textup{\ref{thm:mee}}, 
there exists 
$C>0$ such that
\begin{align*}
&|\mathrm{Im}\,a_{\psi}(x,\xi)-h^{-1-1/s}H_p \psi(x,\xi)| \\
&\le 
(h^{-2/s} \mu^2 +h^{-1}\mu^{\sigma-1/s+1}+\mu^{\sigma-1/s}) 
\quad \text{on} \quad \mathrm{supp}\,[f]
\end{align*}
for $h, \mu \in (0,1]$ with 
$h/\mu \le d.$
\end{lemma}
We can show Lemma \ref{lem:taylor} by a Taylor expansion of 
$\widetilde{a}(x+z_0, \xi+\zeta_0)\,(z_0 \in S_{h^{1-1/s}\nu}, 
\zeta_0 \in \mathbb C^n).$ 
We omit the proof. 
Combining Theorem \ref{thm:mee} and Lemma \ref{lem:taylor}
with the relation $\sigma \le 1/s$, 
we have the following estimate.
\begin{corollary}\label{cor:mee-i}
Under the same assumptions as in Theorem $\ref{thm:mee}$, 
there exists 
$C>0$ such that
\begin{align*}
&|\mathrm{Im}\,
\langle e^{\psi/h^{1/s}}Tu, 
f e^{\psi/h^{1/s}} TPu \rangle
+
\langle e^{\psi/h^{1/s}}Tu, 
f(h^{-1-1/s}H_p \psi) e^{\psi/h^{1/s}} Tu \rangle| \\
& 
\le C
(h^{-2/s} \mu^2 +h^{-1}\mu^{\sigma-1/s+1}+\mu^{\sigma-1/s})
\| \sqrt{f} e^{\psi/h^{1/s}} Tu\|^2 \\
&\quad{}+C
(h^{-1} \mu +\mu^{\sigma}+h\mu^{\sigma-1}) \|u\|^2
\end{align*}
for $u \in L^2(\mathbb R^n)$, $h, \mu \in (0,1]$ with 
$h/\mu \le d.$
\end{corollary}

\section{Proof of Theorem \ref{thm:phwf}}\label{sec:thm}

In this section we prove Theorem \ref{thm:phwf}. 
Using microlocal exponential estimates obtained 
in the previous section, 
we can prove it in the same way as the analytic case \cite{MNS06}.
We use the notation 
\[
 B(x_0,\xi_0 ; a, b)=\{
(x, \xi) \in \mathbb R^{2n}
:
|x-x_0|<a, |\xi-\xi_0|<b
\}.
\]

\vspace*{12pt}

\noindent{\it Proof of Theorem $\ref{thm:phwf}$}
\enskip
We suppose 
$\gamma=\{(y(t),\eta(t)): t \in \mathbb R\}$
is backward nontrapping and 
denote the asymptotic momentum by 
$\eta_-=\lim_{t \to -\infty} \eta(t).$

It follows from \eqref{eqn:phwf-1} that 
there exist a conic neighborhood $\Gamma$ of 
$(-t_0 \eta_-, \eta_-)$ 
and $\delta_0>0$ satisfying
\[
\| e^{ \delta_0 (|x|^{1/s}+|\xi|^{1/s})} T_{1,1} u_0 
\|_{L^2(\Gamma)}<+\infty.
\]
Then we can find a sufficiently small $\delta\in (0,|\eta_-|)$
and $C>0$ such that 
\[
 \|T_{1,1}u_0\|_{L^2( 
B(-h^{-1}t_0 \eta_-,h^{-1}\eta_-; h^{-1}t_0 \delta, h^{-1}\delta)
)} \le Ce^{-\delta/h^{1/s}}
\quad (0< h \le 1),
\]
which is equivalent to 
\begin{equation}\label{eqn:4-1}
\|T_{h, h} u_0\|_{L^2(
B(-h^{-1}t_0 \eta_-, \eta_-; h^{-1}t_0 \delta, \delta)
)} \le Ce^{-\delta/h^{1/s}}
\quad (0< h \le 1).
\end{equation}
Pick up
$\chi_1 (r) \in C^{\infty}(\mathbb R)$ 
satisfying
\[
\chi_1 (r)=
\begin{cases}
1 & (|r| \le 1/2) \\
0 & (|r| \ge 1)
\end{cases},
\quad 
r \chi_1'(r) \le 0.
\]
We set $\delta_1=\delta/4.$ 
Following \cite{MNS06}, we define the weight function by 
\begin{align*}
\psi(t, x, \xi)&=\delta \varphi(t/h,x,\xi), \\ 
\varphi(t, x, \xi)&=\chi_1 \left(\frac{|x-t\xi|}{\delta_1 |t|}\right)
\chi_1 \left(\frac{|\xi-\eta_-|}{\delta_1}\right)
\end{align*}
for $t<0$. 
We recall that 
\begin{align}\label{eqn:4-2}
&\partial_t \varphi+H_p \varphi \le C|t|^{-1-\sigma}, \\
&\mathrm{supp}\,[\varphi(t, \cdot, \cdot)]
\subset B \left(
t\eta_-, \eta_- ; \frac{\delta|t|}{2}, \frac{\delta}{2}
\right). \label{eqn:4-3} 
\end{align}
Combining \eqref{eqn:4-1}, \eqref{eqn:4-3} and the definition of 
$\psi$, we have 
\begin{equation}\label{eqn:4-4}
\|
e^{\psi(t, \cdot, \cdot)/h^{1/s}}T_{h, h}u_0
\| \le C <+\infty
\quad
(|t+t_0| \ll 1). 
\end{equation}
We next set 
\begin{equation}
\mu(t)=-\frac{t_0}{t}h \quad (t<0). 
\end{equation}
It is easy to see that $\mu'(t)>0$ and 
$\mu(t)\ge h$ for $t \in [-t_0, 0)$. 
Moreover, we have 
\begin{align*}
&|\partial_x^{\alpha}\psi|=\mathscr{O}(\mu^{|\alpha|}), \quad 
|\partial_{\xi}^{\alpha}\psi|
=\mathscr{O}(1), \\  
& \mathrm{supp}\,[\psi(t,\cdot,\cdot)]
\subset \{
(x,\xi)\in \mathbb R^{2n}:
C^{-1} \le |\mu x| \le C, C^{-1} \le |\xi| \le C
\}
\end{align*}
for some $C>1$. 
Especially, choosing $\delta$ small again if necessary, 
we see that $\psi(t,x,\xi)$ satisfies (W1). 

Let 
$\chi_2(r) \in C^{\infty}(\mathbb R)$
be a function satisfying
\begin{align*}
\chi_2(r)=
\begin{cases}
1 \quad (A^{-1} \le r \le A), \\
0 \quad (r \le (2A)^{-1}, (2A) \le r),
\end{cases}
\quad 
0 \le \chi_2(r) \le 1 
\end{align*}
for some $A>0$. 
We set 
\[
f(t, x, \xi)=\chi_2(|\mu(t)x|) \chi_2(|\xi|) 
\]
with sufficiently large $A$ so that $f \equiv 1$ on 
$\mathrm{supp}\,[\psi].$ 
We can verify that $f(t,x,\xi)$ satisfies (W2).

For 
$t \in [-t_0, 0)$, 
we set
\[
 F(t)=
\|
\sqrt{f(t,\cdot,\cdot)}e^{\psi(t,\cdot,\cdot)/h^{1/s}}
T_{h, \mu(t)} u(t+t_0)
\|^2.
\]
Then $F(t)$ satisfies the following differential inequality. 
\begin{lemma}\label{lem:di}
There exists $C>0$ such that 
\begin{equation}\label{eqn:di-1}
\frac{d}{dt}F(t) \le
A(t)F(t)+B(t)\|u(t+t_0)\|^2, 
\end{equation} 
where 
\begin{align*}
&A(t)=C(h^{-1/s+\sigma}|t|^{-1-\sigma}+h^{2-2/s}|t|^{-2}),
\quad
B(t)=C(|t|^{-1}+ h^{2} |t|^{-2})
\end{align*} 
$(t<0, 0<h \le 1).$
\end{lemma}
\noindent{\it Proof of Lemma $\ref{lem:di}$.}\enskip
We write $v(t)=u(t+t_0)$ and $T=T_{h, \mu(t)}.$
By a simple computation, we have
\begin{align*}
&\frac{d}{dt}F(t) \\
&=
\langle e^{\psi/h^{1/s}}T (-iPv), fe^{\psi/h^{1/s}}Tv \rangle
+\langle e^{\psi/h^{1/s}}Tv, fe^{\psi/h^{1/s}}T(-iP)v \rangle \\
&\quad{}+\left\langle e^{\psi/h^{1/s}}Tv, 
2h^{-1/s} f \frac{\partial \psi}{\partial t}
e^{\psi/h^{1/s}}Tv \right\rangle \\
&\quad {}+\left\langle e^{\psi/h^{1/s}} 
\left[ \frac{\partial}{\partial t}, T \right]v, 
fe^{\psi/h^{1/s}}Tv \right\rangle
+\left\langle e^{\psi/h^{1/s}}Tv, 
fe^{\psi/h^{1/s}}
\left[ \frac{\partial}{\partial t}, T \right]
v \right\rangle \\
&\quad{}+\left\langle e^{\psi/h^{1/s}}Tv, 
\frac{\partial f}{\partial t}
e^{\psi/h^{1/s}}Tv \right\rangle \\
&=F_1(t)+F_2(t)+F_3(t)+F_4(t)
\end{align*}
We first consider $F_1(t)$ and $F_2(t)$. 
It follows from Corollary \ref{cor:mee-i} that we have 
\begin{align*}
F_1(t)
&=-2\,\mathrm{Im}\,
\langle
e^{\psi/h^{1/s}}Tv, fe^{\psi/h^{1/s}}TPv 
\rangle \\
&=2 
\langle
e^{\psi/h^{1/s}}Tv, f(h^{-1-1/s}H_p \psi) e^{\psi/h^{1/s}}Tv 
\rangle
+r(t),
\end{align*}
where 
\begin{equation}\label{eqn:di-2}
\begin{split}
r(t) &\le
C(h^{2-2/s}|t|^{-2}+h^{\sigma-1/s}|t|^{-\sigma+1/s-1})F(t) \\
&\quad{}+C(|t|^{-1}+h^{\sigma}|t|^{-\sigma})\| v \|^2.
\end{split}
\end{equation}
Therefore, using \eqref{eqn:4-2}, we have 
\begin{equation}\label{eqn:di-3}
\begin{split}
&F_1(t)+F_2(t) \\
&=2h^{-1-1/s}
\left\langle
e^{\psi/h^{1/s}}Tv, 
f\left(h\frac{\partial \psi}{\partial t}+H_p \psi \right) 
e^{\psi/h^{1/s}}Tv 
\right\rangle
+r(t) \\
& \le 
C h^{\sigma-1/s}|t|^{-1-\sigma}F(t)+r(t).
\end{split} 
\end{equation}
We next consider $F_3(t).$ 
Replacing $\psi$ with $h^{1-1/s}\psi$ in the proof of 
\cite[Lemma 4.1]{MNS06}, we can show that 
\begin{equation}\label{eqn:di-4}
\begin{split}
F_3(t)& 
\le C h^{1-2/s}\mu'(t)
\| \sqrt{f} e^{\psi/h^{1/s}} T v \|^2 + C h \mu'(t)\|v \|^2  \\
& \le C' h^{2-2/s}|t|^{-2} F(t)
+C' h^2 |t|^{-2} \| v  \|^2.
\end{split}
\end{equation}
Finally, it is easy to see that 
\begin{equation}\label{eqn:di-5}
F_4(t) \le C |t|^{-1} \| v \|^2.
\end{equation}
Combining 
\eqref{eqn:di-2}, 
\eqref{eqn:di-3}, 
\eqref{eqn:di-4}, 
\eqref{eqn:di-5}, 
we deduce \eqref{eqn:di-1}.
\qed

\vspace*{12pt}

It follows from Gronwall's inequality that we have
\[
 F(t) \le
e^{\int_{-t_0}^t A(\tau)\,d\tau}
\left\{
F(-t_0)+\int_{-t_0}^t B(\tau)d\tau 
\cdot
\sup_{\tau \in [0, t_0]} \| u(\tau) \|^2
\right\} \quad (-t_0 \le t<0).
\]  
This shows that for every $t \in (-t_0, 0)$ 
there exist $C_1>0$, $C_2>0$ such that
\[
F(t)\le C_1 \exp \,( C_2 h^{\sigma-1/s})
\quad (0 < h \le 1). 
\]
Since
\[
\psi(t,x,\xi)=\delta \quad \text{on} \quad  
B \left(
h^{-1}t \eta_-, \eta_-; 
h^{-1}\frac{|t|\delta_1}{4}, \frac{\delta_1}{4} 
\right), 
\] 
we deduce that 
\[
(t\eta_-, \eta_-) \notin \mathrm{HWF}_s (u(t+t_0)) 
\quad (-t_0<t<0), 
\]
which implies \eqref{eqn:phwf-2}.

We next show \eqref{eqn:phwf-3}. 
It is easy to see that for $R \gg 1$ we have 
\begin{align*}
\int_{-t_0}^{-Rh} A(t) \, dt \le \frac{\delta}{2}h^{-1/s}, 
\quad 
\int_{-t_0}^{-Rh} B(t) \, dt \le C|\log h|. 
\end{align*}
Then it follows that
\begin{align*}
F(-Rh)
&=\|
\sqrt{f(-Rh,\cdot,\cdot)}
e^{\psi(-Rh,\cdot,\cdot)/ h^{1/s}}
T_{h, t_0/R} u(t_0-Rh)
\|^2 \\
& \le C e^{\delta/h^{1/s}}
\end{align*}
$(0 < h \le 1)$ for some $C>0.$
It also can be verified that for $R \gg 1$ we have 
\[
\psi(-Rh,x,\xi)=\delta \quad \text{on} \quad
B \left( \gamma(-R); \frac{R\delta_1}{8}, \frac{\delta_1}{8}
\right). 
\]
Then we see that 
\[
\|T_{h, t_0/R}u(t_0-Rh) \|
_{L^2 (B(\gamma(-R);R\delta_1/8, \delta_1/8))}^2 
\le C e^{-\delta/h^{1/s}} \quad (0<h \le 1). 
\]
It is remarked that in the above argument 
we can replace $-t_0$ with $t$ 
in a sufficiently small neighborhood of $-t_0$ 
by virtue of \eqref{eqn:4-4}. 
Then we can find $\varepsilon_0>0$ such that 
\[
\|T_{h, t_0/R} u(\tau-Rh) \|
_{L^2 (B(\gamma(-R);R\delta_1/8, \delta_1/8))}^2 
\le C e^{-\delta/h^{1/s}} 
\]
for $h \in (0,1]$ and 
$\tau \in [-t_0-\varepsilon_0, -t_0+\varepsilon_0].$ 
In particular, we have 
$\gamma(-R) \notin \mathrm{WF}_s (u(\tau))$ 
for
$\tau \in [-t_0-\varepsilon_0/2, -t_0+\varepsilon_0/2].$ 
By the propagation theorem of the microsupport in Gevrey classes, 
which can be shown by the same way as the analytic case
\cite[Lemma 4.3]{MNS06} 
with the aid of Theorem \ref{thm:mee}, 
we can find 
$\varepsilon_1>0$
such that for any 
$t \in \mathbb R$
and 
$\tau \in [-t_0-\varepsilon_1, -t_0+\varepsilon_1]$
there exist 
$C>0$, $\delta'>0$, $\varepsilon>0$ satisfying
\[
\|T_{h, t_0/R} u(\tau-Rh+ht) \|
_{L^2 (B(\gamma(t-R);\delta', \delta'))} 
\le C e^{-\varepsilon /h^{1/s}} 
\quad (0 < h \le 1).
\]
In particular, we have $\gamma(t) \notin \mathrm{WF}_s (u(\tau))$
for any $t \in \mathbb R$ and any 
$\tau \in [-t_0-\varepsilon_1/2, -t_0+\varepsilon_1/2].$
This shows \eqref{eqn:phwf-3}, 
and completes the proof.
\qed

\section{Proof of Corollary \ref{cor:mix}}\label{sec:cor}

This section deals with a microlocal smoothing property
for the initial data
with mixed momentum condition.  
We prove the following lemma 
which immediately implies Corollary \ref{cor:mix}.
\begin{lemma}\label{lem:mix}
Assume that \textup{(A)} holds 
and that $\gamma$ is backward nontrapping. 
Let $\eta_-$ be the asymptotic momentum 
as $t$ tends to $-\infty$.
Assume that there exist 
$\psi(\xi) 
\in C^{\infty}(\mathbb R^n)$
which equals to $1$ in a conic neighborhood of $\eta_-$ and  
$\varepsilon_0>0$, $A_0>0$, $A_1>0$ satisfying
\begin{equation}\label{eqn:mix-1}
\|
(x \cdot D_x)^l \psi(D_x) u_0
\|_{L^2 (\Gamma_{\varepsilon_0})}
\le A_0 A_1^l l!^{2s}
\quad (l \in \mathbb N). 
\end{equation} 
Then we have 
$(-t_0 \eta_-, \eta_-)\notin \mathrm{HWF}_s (u_0)$
for any $t_0>0$.  
\end{lemma}

We give the proof of Lemma \ref{lem:mix} 
following \cite{MRZ}. 
We denote by $B_r (X_0)$ 
the open ball in $\mathbb R^d\,(d=n \,\text{or}\, 2n)$
of radius $r>0$ with centered at $X_0$, that is, 
\[
 B_r (X_0)=\{X \in \mathbb R^d : |X-X_0| <r\}.
\] 
Our goal is to show that for an arbitrarily fixed 
$t_0>0$
there exist 
$C>0$, 
$\delta>0$, 
$r_0>0$ 
and 
$\lambda_0>0$
satisfying 
\begin{equation}\label{eqn:5-1}
|T_{1,1}u_0 (x, \xi)|
\le C e^{-\delta \lambda^{1/s}}
\quad \text{for} \quad 
\lambda \ge \lambda_0, \,
(x, \xi) \in B_{\lambda r_0}( (-\lambda t_0 \eta_-, \lambda\eta_-))
\end{equation}
under the assumption \eqref{eqn:mix-1}.
To show this, we introduce the following operators. 
For 
$\lambda>0$, we set 
\begin{align*}
&S u_0(x,\xi;\lambda)=T_{1,1}\left[
\chi \left(\frac{x}{\lambda}+t_0 \eta_- \right)
\varphi \left(\frac{D_x}{\lambda}-\eta_-\right)
u_0
\right], \\
&\widetilde{S} u_0(x,\xi;\lambda)
=T_{1,1}
\left[
\chi \left(\frac{x}{\lambda}+t_0 \eta_- \right)
\varphi \left(\frac{D_x}{\lambda}-\eta_-\right)
w_0
\right], \\
&w_0=\chi_0 \left(\frac{x}{\lambda}+t_0 \eta_- \right)v_0, 
\quad
v_0=\psi(D_x)u_0,
\end{align*} 
where $\chi$, $\chi_0$, $\varphi \in C^{\infty}(\mathbb R^n)$
satisfy
\begin{align*}
&|\partial_x^{\alpha}\chi(x)| \le A_2^{1+|\alpha|}\alpha!^s,
\quad 0 \le \chi(x) \le 1 
\quad
(\alpha \in \mathbb Z_+^n, x \in \mathbb R^n),
\\  
&\chi_0, \varphi \quad \text{also satisfy the same estimates}, \\
&\chi(x)=
\begin{cases}
1 \quad (|x| \le \varepsilon_1) \\
0 \quad (|x| \ge 2\varepsilon_1)
\end{cases}, 
\quad
\chi_0(x)=
\begin{cases}
1 \quad (|x| \le 2\varepsilon_1) \\
0 \quad (|x| \ge 3\varepsilon_1)
\end{cases}, \\
&
\varphi(\xi)=
\begin{cases}
1 \quad (|\xi| \le \varepsilon_2) \\
0 \quad (|\xi| \ge 2\varepsilon_2)
\end{cases},
\end{align*}
and 
$\varepsilon_1>0$, $\varepsilon_2>0$ 
are sufficiently small so that 
\begin{align} \label{eqn:5-2}
&\psi =1 \quad \text{on} \quad \mathrm{supp}\,
\left[\varphi \left(
\frac{(\cdot)}{\lambda}-\eta_-
\right)
\right] \quad (\lambda>0), \\
&
\left|\frac{\lambda^2}{x \cdot \xi}\right|
\le A_3 \quad \text{for} \quad 
\lambda>0, \,
(x,\xi) \in 
B(-\lambda t_0 \eta_-, \lambda \eta_-;
2\varepsilon_1 \lambda, 2 \varepsilon_2 \lambda ), 
\label{eqn:5-3} \\
&
\mathrm{supp}\,\left[
\chi_0 \left(
\frac{(\cdot)}{\lambda}+t_0 \eta_-
\right)
\right]
\subset \Gamma_{\varepsilon_0}
\quad \text{for} \quad \lambda \gg 1.
\label{eqn:5-4}
\end{align} 
The proof of Lemma \ref{lem:mix}
is divided into three lemmas.
\begin{lemma}\label{lem:mix1}
There exist $C>0, \delta>0$ such that 
\begin{equation} \label{eqn:mix1-1}
|T_{1,1}u_0 (x,\xi)-Su_0(x,\xi; \lambda)|
\le 
C e^{-\delta \lambda^2}\|u_0 \|
\end{equation} 
for $\lambda>0$, $(x,\xi) 
\in 
B
\left(
-\lambda t_0 \eta_-, \lambda \eta_-;
\lambda \varepsilon_1/2, 
\lambda \varepsilon_2/2 \right).$
\end{lemma}
\begin{lemma}\label{lem:mix2}
There exist $C>0, \delta>0$ such that 
\begin{equation} \label{eqn:mix2-1}
|S u_0 (x,\xi; \lambda)-\widetilde{S} u_0(x,\xi; \lambda)|
\le 
C e^{-\delta \lambda^{1/s}}\|u_0 \|
\end{equation} 
for $\lambda>0$, $x,\xi \in \mathbb R^n$.
\end{lemma}
\begin{lemma}\label{lem:mix3}
Assume \eqref{eqn:mix-1}. Then
there exist $C>0$, $\delta>0$ and $\lambda_0>0$ 
such that 
\begin{equation} 
\label{eqn:mix3-1}
|\widetilde{S} u_0(x,\xi; \lambda)|
\le 
C e^{-\delta \lambda^{1/s}}
\end{equation} 
for $\lambda \ge \lambda_0$, $x,\xi \in \mathbb R^n$.
\end{lemma}
\noindent{\it Proof of Lemma $\ref{lem:mix}$.}\enskip 
It follows from \eqref{eqn:mix1-1}, 
\eqref{eqn:mix2-1} and \eqref{eqn:mix3-1}
that we obtain \eqref{eqn:5-1}, 
which implies $(-t_0 \eta_-, \eta_-) \notin \mathrm{HWF}_s (u_0).$
\qed

\vspace*{12pt}

In the rest of this section, 
we prove these three lemmas. 
We denote the Fourier transform of 
$u(x) \in \mathscr{S}(\mathbb R^n)$
by 
\[
\widehat{u}(\xi)
=\frac{1}{(2\pi)^{n/2}}
\int_{\mathbb R^n}
e^{-ix\cdot \xi}u(x)\,dx.
\]

\vspace*{12pt}

\noindent{\it Proof of Lemma $\ref{lem:mix1}$.}
\enskip
Since $1=\chi \varphi+(1-\chi)\varphi+1-\varphi$, 
we have
\begin{align*}
&T_{1, 1}u_0 (x,\xi)-S u_0 (x,\xi;\lambda) \\
&=c_{1,1}
\int e^{i(x-y)\cdot \xi-|x-y|^2/2}
\left(
1-\chi \left(\frac{y}{\lambda}+t_0 \eta_- \right)
\right)
\varphi \left(
\frac{D_y}{\lambda}-\eta_- \right)
u_0 (y)\,dy \\
&\quad {}
+c_{1,1}
\int e^{i(x-y)\cdot \xi-|x-y|^2/2}
\left(1- \varphi \left(
\frac{D_y}{\lambda}-\eta_- \right)\right)
u_0 (y)\,dy \\
&=I+II.
\end{align*}
We note that if $|x+\lambda t_0 \eta_-|<\lambda \varepsilon_1/2$, 
then we have $|x-y| \ge \lambda \varepsilon_1/2$
on the support of the integrand of the term $I$.
Therefore, using the Schwarz inequality, we have
\begin{align*}
|I| & \le c_{1,1} e^{-\varepsilon_1^2 \lambda^2 / 16}
\int e^{-|x-y|^2/4}
\left|
(1-\chi) \varphi \left(
\frac{D_y}{\lambda}-\eta_- \right)
u_0 (y)
\right|
\,dy \\
& 
\le  c_{1,1} e^{-\varepsilon_1^2 \lambda^2 / 16}
(2\pi)^{n/4}
\|u_0\|
\end{align*}
We can treat the term $II$ in the same way. Indeed, since
$T_{1,1}u(x,\xi)=e^{ix\cdot \xi}T_{1,1}\widehat{u}(\xi,-x)$, 
we have 
\[
II
=e^{ix \cdot \xi} \cdot 
c_{1,1}
\int e^{i(\xi-\eta) \cdot (-x)-|\xi-\eta|^2/2}
\left(1- \varphi \left(
\frac{\eta}{\lambda}-\eta_- \right)\right)
\widehat{u_0} (\eta)\,d\eta.
\]
Then it follows that 
\[
|II| \le  c_{1,1} e^{-\varepsilon_2^2 \lambda^2 / 16}
(2\pi)^{n/4}
\| \widehat{u_0} \|,
\]
which completes the proof.
\qed

\vspace*{12pt}

\noindent{\it Proof of Lemma $\ref{lem:mix2}$.}
\enskip
We set 
\[
 I v_0(y;\lambda)
=\chi 
\left(
\frac{y}{\lambda}+t_0 \eta_- \right)
\varphi 
\left(
\frac{D_y}{\lambda}-\eta_-
\right)
\left(
\chi_0
\left(\frac{y}{\lambda}+t_0 \eta_- \right)-1
\right)
v_0(y).
\]
Then, using \eqref{eqn:5-2}, we can verify that
\begin{equation}\label{eqn:mix2-2}
\widetilde{S}u_0-S u_0
=T_{1,1}[I v_0].
\end{equation}
It suffices to show that there exist 
$B_0>0$, $B_1>0$ such that
\begin{equation}\label{eqn:mix2-3}
\lambda^{2N}
|I v_0(y,\lambda)|
\le 
B_0 B_1^N N^{2sN} \lambda^{n/2}
\|v_0\| \chi \left(
\frac{y}{\lambda}+t_0 \eta_-
\right)
\end{equation}
for $N \in \mathbb N, \lambda>0, y \in \mathbb R^n.$
Indeed, combining \eqref{eqn:mix2-2}, \eqref{eqn:mix2-3}
and $N^N \le e^N N!$, we have 
\begin{align*}
|\widetilde{S}u_0-S u_0|
& \le B_2 \lambda^{n/2} B_1^N \lambda^{-2N} N^{2sN} \| v_0 \| \\
& \le B_2 \lambda^{n/2} \left\{
\frac{N!}{(B_3 \lambda^{1/s})^N}
\right\}^{2s} \| v_0\|,  
\end{align*}
where $B_3^{-1}=B_1^{1/(2s)}e.$
Choosing $N=[B_3 \lambda^{1/s}]$, we deduce \eqref{eqn:mix2-1}.

We can prove \eqref{eqn:mix2-3} in the same way as 
\cite[Lemma III.4]{MRZ} .
We recall that 
\begin{align*}
Iv_0 (y,\eta)
&=\frac{1}{(2\pi)^n}\chi\left(\frac{y}{\lambda}+t_0 \eta_-\right) \\
&\quad \times 
\int \!\!\! \int e^{i(y-z)\cdot \eta}
\varphi\left(\frac{\eta}{\lambda}-\eta_-\right) 
\left(
\chi_0 \left(\frac{z}{\lambda}+t_0 \eta_-\right)-1
\right)
v_0(z)\,dz d\eta.
\end{align*}
Since $\chi_0 =1$ on $\mathrm{supp}\,[\chi]$, we can replace
$\chi_0 \left(z/\lambda+t_0 \eta_-\right)-1$
in the above with 
\[
\chi_0 \left(\frac{z}{\lambda}+t_0 \eta_-\right)-
\chi_0 \left(\frac{y}{\lambda}+t_0 \eta_-\right),
\]
which equals to 
\begin{align*}
&\sum_{1 \le |\alpha| \le N-1}
\frac{(z-y)^{\alpha}}{\alpha! \lambda^{|\alpha|}}
(\partial_y^{\alpha} \chi_0)\left(\frac{y}{\lambda}+t_0 \eta_-\right)
\\
&\quad{}
+N \sum_{|\alpha|=N} 
\frac{(z-y)^{\alpha}}{\alpha! \lambda^{|\alpha|}}
\int_0^1
(1- \theta)^{N-1}
(\partial_z^{\alpha} \chi_0)
\left(\frac{\theta z+(1-\theta)y}{\lambda}+t_0 \eta_-\right)\,
d\theta.
\end{align*}
Using 
$\chi(y) \cdot (\partial_y^{\alpha}\chi_0)(y)=0 \, (|\alpha| \ne 0)$ 
and 
$(z-y)^{\alpha}e^{i(y-z)\cdot \eta}
=(i\partial_{\eta})^{\alpha}e^{i(y-z)\cdot \eta}$, 
we deduce that 
\begin{equation}
\label{eqn:mix2-4}
\begin{split}
&I v_0 (y;\lambda) \\
&=\frac{1}{(2\pi)^{n}}
\chi \left(\frac{y}{\lambda}+t_0 \eta_-\right)
\sum_{|\alpha|=N}
\frac{N}{\alpha! \lambda^{2|\alpha|}}
\int \!\!\! \int e^{i(y-z)\cdot \eta}
(D_{\eta}^{\alpha} \varphi)\left(\frac{\eta}{\lambda}-\eta_-\right) \\
&\qquad \times
\left\{
\int_0^1
(1- \theta)^{N-1}
(\partial_z^{\alpha} \chi_0)
\left(\frac{\theta z+(1-\theta)y}{\lambda}+t_0 \eta_-\right)\,
d\theta
\right\}
v_0(z)\,dzd\eta 
\end{split}
\end{equation}
Here we use the integration by parts with respect to $\eta$.
By the Schwarz inequality, the Plancherel formula and the properties of 
$\varphi$ and $\chi_0$, we can find $B_4>0$ such that 
\begin{align*}
&\left|
\int \!\!\! \int 
e^{i(y-z) \cdot \eta} 
(D_{\eta}\varphi)\left(\frac{\eta}{\lambda}-\eta_- \right)
\left\{\int_0^1 (
\cdots) \,d\theta \right\} v_0(z) \,dz d\eta
\right| \\
&\le B_4^{1+|\alpha|}\lambda^{n/2}|\alpha|^{2s|\alpha|}
\|v_0\|, 
\end{align*}
which and \eqref{eqn:mix2-4} imply \eqref{eqn:mix2-3}.
\qed

\vspace*{12pt}

Finally we give the proof of Lemma \ref{lem:mix3}.
We set 
\begin{equation}
\label{eqn:5-5}
\begin{split}
I_{J,N,k}u(y; \lambda)
&=\frac{1}{(2\pi)^{n/2}}
\chi \left( \frac{y}{\lambda}+t_0 \eta_-\right) \\
& \quad \times \int e^{iy \cdot \eta}
\left(\frac{\lambda^2}{y \cdot \eta}\right)^N 
\left[
(\eta \cdot \partial_{\eta})^k 
\left\{
\varphi \left( \frac{\eta}{\lambda}-\eta_- \right)
\right\}
\right] (y \cdot \eta)^J
\widehat{u}(\eta)\,d\eta,
\end{split} 
\end{equation}
then we have 
\begin{equation}\label{eqn:5-6}
\lambda^{2N}\widetilde{S}u_0 (x,\xi;\lambda)=T_{1,1}[I_{N,N,0}w_0].
\end{equation}
It is remarked that we can divide by $y \cdot \eta$ in \eqref{eqn:5-5}
since we have \eqref{eqn:5-3}. 
Corresponding to 
\cite[Lemma III.9]{MRZ} 
and 
\cite[Lemma III.8]{MRZ}, 
we have the following results.
\begin{lemma}
\label{lem:mix4}
There exist $D_0>0$, $D_1>0$ and $D_2>0$ such that 
\begin{equation}
\label{eqn:mix4-1}
\| I_{J,N,k}u \| \le D_0 D_1^N D_2^J \lambda^n 
\sum_{l=0}^J N^{2s(J-l)+sk} \| (y\cdot D_y)^l u \| 
\end{equation}
for $J \ge 0, N \ge J, k \le N-J$ and 
$u \in \mathscr{S}(\mathbb R^n).$ 
\end{lemma}
\begin{lemma}\label{lem:mix5}
Assume \eqref{eqn:mix-1}. 
Then there exist $E_0>0$, $E_1>0$ and $\lambda_1>0$ such that 
\begin{equation}
\label{eqn:mix5-1}
\|
(y \cdot D_y)^l w_0
\|
\le 
E_0 E_1^l l^{2sl} 
\end{equation}
for $l \in \mathbb N$, $\lambda \ge \lambda_1.$ 
\end{lemma}
Using these lemmas, 
we can estimate $\widetilde{S}u.$

\vspace*{12pt}

\noindent{\it Proof of Lemma \textup{\ref{lem:mix3}}.}\enskip
Combining Lemma \ref{lem:mix4} with $J=N$, $k=0$ 
and Lemma \ref{lem:mix5}, 
we see that 
\begin{align*}
\|I_{N,N,0}w_0 \| 
&\le 
D_0 D_1^N D_2^N \lambda^n 
\sum_{l=0}^N N^{2s(N-l)} \| (y \cdot D_y)^l w_0 \| \\
&\le
D_0 D_1^N D_2^N \lambda^n 
\sum_{l=0}^N N^{2sN} E_0 E_1^l \\
& \le 
D_3 D_4^N N^{2sN} \lambda^n 
\end{align*}
for $\lambda \ge \lambda_1.$ 
Then, using \eqref{eqn:5-6}, we deduce that 
\[
 |\lambda^{2N} \widetilde{S}u_0(x,\xi;\lambda)|
\le D_5 D_4^N N^{2sN} \lambda^n
\]
for $\lambda \ge \lambda_1$, which implies 
\eqref{eqn:mix3-1}.
\qed

\vspace*{12pt}

The proofs of Lemma \ref{lem:mix4} and Lemma \ref{lem:mix5}
are similar to those of 
\cite[Lemma III.9]{MRZ} and \cite[Lemma III.8]{MRZ} respectively.
We only mention the sketch of them.

\vspace*{12pt}

\noindent{\it Proof of Lemma $\ref{lem:mix4}.$}\enskip
We use an induction argument for $J$. 
First we can verify that there exist $M_0>0$, $M_1>0$ such that 
\[
\left|
(\eta \cdot \partial_{\eta})^k 
\left\{
\varphi \left(\frac{\eta}{\lambda}-\eta_-\right)
\right\}
\right|
\le M_0 M_1^k k^{sk} \quad (k \in \mathbb N).
\]
Combining this and \eqref{eqn:5-3}, 
we obtain \eqref{eqn:mix4-1} with $J=0.$

For the case of general $J$, we make use of 
\begin{align*}
I_{J,N,k}u&=in I_{J-1, N,k}u-iN I_{J-1,N,k}u+iI_{J-1,N,k+1}u \\
&\quad{}+i(J-1)I_{J-1,N,k}u-iI_{J-1,N,k}(nu+y \cdot \partial_y u).
\end{align*}
We omit the details.
\qed

\vspace*{12pt}

\noindent{\it Proof of Lemma \textup{\ref{lem:mix5}}.}\enskip
We recall that
\begin{equation}
\label{eqn:mix5-2}
(y\cdot D_y)^l w_0
=\sum_{k=0}^l \binom{l}{k} (y \cdot D_y)^k 
\left\{
\chi_0 \left(\frac{y}{\lambda}+t_0 \eta_-\right)
\right\}
(y\cdot D_y)^{l-k} v_0.
\end{equation}
It can be verified that there exist $M_2>0$, $M_3>0$ such that 
\begin{equation}
\label{eqn:mix5-3}
\left|
(y\cdot D_y)^k
\left\{
\chi_0 \left(\frac{y}{\lambda}+t_0 \eta_-\right)
\right\}
\right|
\le M_2 M_3^k k^{sk} 
\quad (k \in \mathbb N).
\end{equation}
On the other hand, using \eqref{eqn:mix-1} and \eqref{eqn:5-4}, 
we have
\begin{equation}
\label{eqn:mix5-4}
\| (y \cdot D_y)^l v_0 \|
_{ L^2( B_{3 \varepsilon_1 \lambda}(-\lambda t_0 \eta_-))}
\le A_0 A_1^l l!^{2s} 
\end{equation}
for $l \in \mathbb N$ and $\lambda \gg 1.$
Combining 
\eqref{eqn:mix5-2},
\eqref{eqn:mix5-3} 
and 
\eqref{eqn:mix5-4},
we deduce 
\eqref{eqn:mix5-1}.
\qed

\end{document}